\documentclass{amsart}
\usepackage{amscd,amssymb,hyperref}
\usepackage[all]{xy}
\theoremstyle{plain}
\newtheorem{prop}[subsection]{Proposition}
\newtheorem{thm}[subsection]{Theorem}
\newtheorem{lem}[subsection]{Lemma}
\newtheorem{cor}[subsection]{Corollary}
\newtheorem{remark}[subsection]{\bfseries\itshape Remark}

\theoremstyle{definition}
\newtheorem{defn}[subsection]{Definition}
\newtheorem{exm}[subsection]{Example}
\numberwithin{equation}{section}

\begin{document}

\title{Braid groups, free groups, and the loop space of the $2$-sphere}
\author[F.R.~Cohen]{F.R. Cohen$^{*}$}
\address{Department of Mathematics, University of Rochester,
Rochester, New York 14627 } \email{cohf@math.rochester.edu}
\urladdr{www.math.rochester.edu/\~{}cohf}

\author[J.~Wu]{J.~Wu$^\dag$}
\address{Department of Mathematics, National University of Singapore,
Singapore} \email{matwuj@nus.edu.sg}
%%\urladdr{}
\thanks{$^{*}$Partially supported by the NSF}
\thanks{$^\dag$Partially supported by a grant
from the National University of Singapore}

%%\subjclass[2000]{Primary~20F36, 55R50; Secondary }
%%20F36 Special aspects of infinite or finite groups
  %%     Braid groups; Artin groups
%%52C35 Discrete geometry
  %%    Arrangements of points, flats, hyperplanes
%%55R50 Fiber spaces and bundles
  %%     Stable classes of vector space bundles, $K$-theory
    %%   Simplicial groups, homotopy theory

\date{\today}


\begin{abstract}
The purpose of this article is to describe connections between the
loop space of the $2$-sphere, Artin's braid groups, a choice of
simplicial group whose homotopy groups are given by modules called
$\mathrm{Lie}(n)$, as well as work of Milnor \cite{ Milnor}, and
Habegger-Lin \cite{Habegger-Lin, Lin} on "homotopy string links".
The current article exploits Lie algebras associated to Vassiliev
invariants in work of T. Kohno \cite{K,K1}, and provides
connections between these various topics.

Two consequences are as follows:
\begin{enumerate}
    \item the homotopy groups of spheres are identified as
    ``natural" sub-quotients of free products of pure braid groups, and
    \item an axiomatization of certain simplicial groups
    arising from braid groups is shown to characterize
    the homotopy types of connected $CW$-complexes.
\end{enumerate}
\end{abstract}

 \maketitle
\section{A tale of two groups plus one more}

In $1924$ E.~Artin \cite{Artin, A} defined the $n$-th braid group
$B_n$ together with the $n$-th pure braid group $P_n$, the kernel
of the natural map of $B_n$ to $\Sigma_n$ the symmetric group on
$n$-letters. It is the purpose of this article to derive
additional connections of these groups to homotopy theory, as well
as some overlaps with algebraic, and topological properties of
braid groups.

This article gives certain relationships between free groups on
$n$ generators $F_n$, and braid groups which serve as a bridge
between different structures. These connections, at the interface
of homotopy groups of spheres, braids, knots, and links, and
homotopy links, admit a common thread given by a simplicial group.

Recall that a simplicial group $\Gamma_{*}$ is a collection of
groups $$\Gamma_0,\Gamma_1, \cdots, \Gamma_n, \cdots $$ together
with face operations $$d_i:\Gamma_n \to\ \Gamma_{n-1},$$  and
degeneracy operations $$s_i:\Gamma_n \to\ \Gamma_{n+1},$$ for $ 0
\leq i \leq n$. These homomorphisms are required to satisfy the
standard simplicial identities.

One example is Milnor's free group construction $F[K]$ for a
pointed simplicial set $K$ with base-point $*$ in degree zero. The
simplicial group $F[K]$ in degree $n$ is the free group generated
by the $n$ simplices $K_n$ modulo the single relation that
$s_0^n(*) = 1$.

In case $K$ is reduced, that is $K$ consists of a single point in
degree zero, the geometric realization of $F[K]$ is homotopy
equivalent to $\Omega \Sigma|K|$ \cite{Milnor}. The first theorem
below addresses one property concerning the simplicial group given
by $F[\Delta[1]]$ where $\Delta[1]$ is the simplicial $1$-simplex,
and a simplicial group given in terms of the pure braid groups
described next.

A second example is given by the simplicial group which in degree
$n$ is given by $\Gamma_n = P_{n+1}$, the $(n+1)$-st pure braid
group, and which is elucidated in \cite{CW,BCWW}. The face
operations are given by deletion of a strand, while the
degeneracies are gotten by ``doubling" of a strand. This
simplicial group is denoted $\mathrm{AP}_*$.

\begin{thm} \label{thm:Looping braids}
The loop space ( as simplicial groups ) of the simplicial group
$\mathrm{AP}_*$, $\Omega(\mathrm{AP}_*)$, is isomorphic to
$F[\Delta[1]]$ as a simplicial group.
\end{thm}

A connection between Artin's braid group and the loop space of the
$2$-sphere is given next using the feature that the second pure
braid group is isomorphic to the integers with a choice of
generator denoted $A_{1,2}$. Notice that the simplicial circle
$S^1$ has a single non-degenerate point in degree $1$ given by
$\langle 0,1\rangle$. Thus there exists a unique map of simplicial
groups
$$\Theta: F[S^1] \to\ \mathrm{AP}_*$$ such that $\Theta(\langle 0,1\rangle) =  A_{1,2}$.
One of the theorems stated in \cite{CW} is as follows.

\begin{thm}\label{thm:Embedding braids}
The morphism of simplicial groups
$$\Theta: F[S^1] \to\ \mathrm{AP}_*$$ is an embedding.
Hence the homotopy groups of $F[S^1]$ are natural sub-quotients of
$\mathrm{AP}_*$, and the geometric realization of quotient
simplicial set $\mathrm{AP}_*/F[S^1]$ is homotopy equivalent to
the $2$-sphere. Furthermore, the image of $\Theta$ is the smallest
simplicial subgroup of $\mathrm{AP}_*$ which contains $A_{1,2}$.
\end{thm}

This theorem gives that the homotopy groups of $F[S^1]$, those of
the loop space of the $2$-sphere, are given as ``natural"
sub-quotients of the braid groups, a result related to work of the
second author \cite{W}. The proof of the above theorem sketched in
\cite{CW} relies heavily on the structure of a Lie algebra arising
from the ``infinitesimal braid relations" as Vassiliev invariants
of braids by work of T.~Kohno \cite{K,K1,K2}, Falk, and Randell
\cite{FR}, as well as work of V.~Drinfel'd \cite{D} on the $KZ$
equations. The precise details in section $9$ here depend heavily
on the specific structure of this Lie algebra.

An example is listed next. The commutator of the braids $x_1$, and
$x_2$ in the third braid group as listed below in section
\ref{sec:On Theta n} represents the Hopf map $\eta:S^3 \to S^2$.
The braid closure of this commutator gives the Borromean rings.

An analogue for all spheres arises at once by taking coproducts of
simplicial groups $\mathrm{AP}_* \vee \mathrm{AP}_*$ which in
degree $n$ is given by the free product $P_{n+1} \amalg  P_{n+1}$.
\begin{cor}\label{cor:all spheres}
The smallest simplicial subgroup of $\mathrm{AP}_* \vee
\mathrm{AP}_*$ which contains $P_{2} \amalg P_{2}$ in degree $1$
is isomorphic to $F[S^1]\vee F[S^1]$. Hence $\Omega S^n$ is a
retract, up to homotopy, of the geometric realization of the
simplicial subgroup $F[S^1]\vee F[S^1]$ for any $n \geq 2$ by the
Hilton-Milnor theorem.
\end{cor}

In addition, the Lie algebraic methods used to prove this theorem
suggest that the methods might be useful to study whether related
maps are faithful.  A ``Lie algebraic/homological" criterion for
testing whether a representation of $P_n$ is faithful is given in
\cite{CP}.

The results above suggest an axiomatization of certain families of
simplicial groups. One application of Theorems \ref{thm:Looping
braids}, and \ref{thm:Embedding braids} is listed next. Let
$\mathcal B$ denote the smallest full sub-category of the category
of reduced simplicial groups which satisfies the following
properties:
\begin{enumerate}
    \item The simplicial group $\mathrm{AP}_*$ is in $\mathcal B$.
    \item If $\Pi$, and $\Gamma$ are in $\mathcal B$, then the
    coproduct $\Pi \vee \Gamma$ is in $\mathcal B$.
    \item If $\Pi$ is in $\mathcal B$, and $\Gamma$ is a
    simplicial subgroup of $\Pi$, then $\Gamma$ is in $\mathcal B$.
     \item If $\Pi$ is in $\mathcal B$, and $\Gamma$ is a
     simplicial
    quotient of $\Pi$, then $\Gamma$ is in $\mathcal B$.
\end{enumerate} To be precise, the authors are unaware of a
specific reference for the definition of the object given by a
simplicial subgroup. There are two natural, and equivalent
definitions given in section \ref{sec:axioms}.

\begin{thm}\label{thm:theorem one}
Let $X(i)$, $i = 1,2$ denote path-connected $CW$-complexes with a
continuous function $f:X(1) \to X(2)$. Then there exist elements
$\Gamma_{X(i)}$, together with a morphism $\gamma:\Gamma_{X(1)}
\to \Gamma_{X(2)}$ in $\mathcal B$ such that the loop space
$\Omega(X(i))$ is homotopy equivalent to the geometric realization
of $\Gamma_{X(i)}$, and the induced map $|\gamma|:|\Gamma_{X(1)}|
\to |\Gamma_{X(2)}|$ is homotopic to $\Omega(f)$.
\end{thm}

One property concerning these connections is that the Lie algebra
obtained from the descending central series of a discrete group
has several disparate features here. One is that the structure of
the resulting Lie algebra gives the method for proving that map in
the proof of Theorem \ref{thm:Embedding braids} is an injection.
The second feature is that this Lie algebra is used to give the
axiomatization given in Theorem \ref{thm:theorem one}. As remarked
above, this Lie algebra is used to characterize Vassiliev
invariants of pure braids \cite{K,K1}. Lastly, this Lie algebra as
well as the Lie algebra obtained from the mod-$p$ descending
central series for a discrete group gives the classical the
Bousfield-Kan spectral sequence, as well as the classical unstable
Adams spectral sequence for which the descending central series is
replaced by the mod-$p$ descending central series.

It is natural to consider other quotients of the simplicial
groups, and Lie algebras here such as those arising in work of
Habegger-Lin, \cite{Habegger-Lin}. That is the purpose of the next
few remarks with one goal given by a construction of natural
quotients of the pure braid groups in Proposition \ref{prop:
descent of the braid group}, as well as \ref{cor: the pullback for
the group Pn}.

A naturality argument provides a bridge between the above
theorems, certain quotients of the braid groups as well as the
modules $\mathrm{Lie}(n)$ related to the cohomology of the pure
braid groups. To make this naturality argument precise, a setting
to provide natural quotients of the braid group is required. This
setting provides an identification of a classical representation
of the pure braid group in the automorphism group of a free group
due to E.~Artin \cite{Artin,A} in terms of a ``universal"
semi-direct product, the so-called holomorph of a group.  The
utility of this information is that it provides a direct method
for constructing natural quotients of the braid groups which fit
within the contexts of homotopy theory, links, and ``homotopy
string links", as well as intermediate analogues. In addition,
features of the holomorph give a direct proof of Artin's classic
result that his representation of the braid group in the
automorphism of a free group is faithful. This proof then extends
at once to other natural quotients. An extension of this
information is applied to certain choices of simplicial groups one
of which arises from ``homotopy string links".

Related results are given in \cite{CW,BCWW}.

\begin{enumerate}
    \item [] {\bf TABLE OF CONTENTS}
\end{enumerate}

\begin{enumerate}
  \item Introduction
  \item On the holomorph of a group, and the pure braid groups as group extensions
  \item Quotients of the pure braid group, and the maps $B_n \to\
Aut(K_n)$
  \item On link homotopy and related homotopy groups
  \item On looping $\mathrm{AP}_*$
  \item On embeddings of residually nilpotent groups
  \item The simplicial structure for $\mathrm{AP}_*$
  \item The Lie algebra associated to the descending central series
for $P_{n+1}$
  \item On $\Theta_n: F_{n} \to\ P_{n+1}$
  \item On the proof of Theorem \ref{thm: embeddings of Lie algebras}
  \item On Vassiliev invariants, the mod-$p$ descending central
series, and the Bousfield-Kan spectral sequence
  \item On braid groups, and axioms for connected $CW$-complexes
  \item Proof of Theorem \ref{thm:theorem one}
  \item Appendix: a sample computation
\end{enumerate}

The first author would like to thank the University of Tokyo, and
the National University of Singapore for support during some of
the work on this paper. The second author would like to thank the
Universit\'e Montpellier II for support during some of the work on
this project.

\section{ On the holomorph of a group, and the pure braid groups as group extensions}
\label{sec:pure braid groups as extensions}

The purpose of this section is to give methods for direct
constructions of quotients of the pure braid groups via natural
universal split group extensions. One application below is to
extend work above on $\Theta_n:F_n \to P_{n+1}$ of Theorem
\ref{thm:Embedding braids} to the setting of homotopy string
links.

Consider a discrete group $\pi$ together with the universal
semi-direct product $\mathrm{Hol}(\pi)$, as elucidated in work of
M.~Voloshina \cite{Voloshina} with some details given below. The
group $\mathrm{Hol}(\pi)$ is ``the natural" split extension of
$\mathrm{Aut}(\pi)$ by $\pi$, $$1 \to\ \pi \to\ \mathrm{Hol}(\pi)
\to\ \mathrm{Aut}(\pi) \to\ 1$$ where $\mathrm{Aut}(\pi)$ is the
automorphism group of $\pi$. More precisely, the group
$\mathrm{Hol}(\pi)$, as a set, is a product $\mathrm{Aut}(\pi)
\times \pi$, but the group structure is defined by the product
$$(f,x)\cdot (g,y) = (f\cdot g, g^{-1}(x) \cdot y)$$ for $f,g$ in
$\mathrm{Aut}(\pi)$, and $x,y$ in $\pi$. Hence, $$(f,1)^{-1}\cdot
(1,y)\cdot (f,1)= (1, f^{-1}(y)).$$

The topological analogue of this construction is the universal
bundle with section having fibre a given $K(\pi,1)$:
$$E\mathrm{Aut}(\pi) \times_{\pi} K(\pi,1) \to B\mathrm{Aut}(\pi) $$ where
$\mathrm{Aut}(\pi)$ acts by ``left translation" on $\pi$. This
construction will be used below to describe groups analogous to
the braid groups obtained from a classical representation due to
Artin.

Let $F_n$ denote the free group on $n$-letters with generators
$x_1,x_2,\cdots, x_n$. Artin gave a homomorphism from the $n$-th
braid group to $\mathrm{Aut}(F_n)$ \cite{A,Artin,B,MKS} $$A:B_n
\to\ \mathrm{Aut}(F_n)$$ which is described below. One consequence
is that the restriction of $A$ to $P_n$ identifies Artin's
representation of $P_{n+1}$ inductively arising from a pull-back
in Corollary \ref{cor: the pullback for the group Pn} of the two
maps $A:P_n \to\ \mathrm{Aut}(F_n)$, and $\mathrm{Hol}(F_n) \to\
\mathrm{Aut} (F_{n})$.

To describe this structure, recall the following results due to
E.~Artin  \cite{Artin,A,B,MKS} where the commutator is defined by
the equation $[a,b] =  a^{-1} \cdot b^{-1} \cdot a \cdot b$. The
braid group on $n$-strands $B_n$ is generated by elements,
$\sigma_i$ for $1 \leq i \leq n-1$ with the well-known relations
\begin{enumerate}
  \item $\sigma_i \sigma_j = \sigma_j \sigma_i$ if $|i-j| \geq 2$,
  and
  \item $\sigma_i \sigma_{i+1} \sigma_i  = \sigma_{i+1} \sigma_{i} \sigma_{i+1}$
  for all $i$.
\end{enumerate} Artin's representation of $B_n$ is stated in the next theorem.

\begin{thm}\label{theorem:Artin's representation}
There is a homomorhism
$$A: B_n \to\ \mathrm{Aut}(F_n)$$ obtained by setting
\begin{enumerate}
    \item $A(\sigma_s)(x_j) = x_j$ if $j \neq s, s+1$,
    \item $A(\sigma_s)(x_s) = x_{s+1}$, and
    \item $A(\sigma_s)(x_{s+1}) = x_{s+1}^{-1}x_sx_{s+1} = x_s[x_{s},x_{s+1}]
    $.
\end{enumerate} Thus

\begin{enumerate}
    \item $A(\sigma_s)^{-1}(x_j) = x_j$ if $j \neq s, s+1$,
    \item $A(\sigma_s)^{-1}(x_s) = x_{s}x_{s+1}x_{s}^{-1} $, and
    \item $A(\sigma_s)^{-1}(x_{s+1}) = x_s$.
\end{enumerate}

\end{thm}

The pure braid group $P_n$ is generated by elements $A_{r,s}$ for
$1 \leq r < s \leq n$. The element $A_{r,s}$ may be thought of as
linking the $s$-th strand around the $r$-th strand for $1 \leq r <
s \leq n$ with a fixed orientation together with an explicit
formula given by
$$A_{r,s} =  \alpha(r,s)\cdot \sigma_r^2\cdot  \alpha(r,s)^{-1}  $$
where $$\alpha(r,s) = \sigma_{s-1}\cdot \sigma_{s-2}\cdot
\sigma_{s-3}\cdots \sigma_{r+1}.$$ Artin gave a complete set of
relations for the pure braid group $P_n$ \cite{A, Artin}, or
\cite{MKS}. A restatement of Artin's relations in terms of
commutators is listed next.

\begin{thm}\label{thm: action of pure braids on a free group}
The group $P_n$ is generated by elements
$$A_{r,s}$$ for $1 \leq r < s \leq n$. A complete set of
relations is given as follows:
\begin{enumerate}

    \item If either $ s<i$, or $ k<r$,
    $A_{r,s}A_{i,k}A_{r,s}^{-1} = A_{i,k}$.

 \item  If $ i<k <s$,
 $A_{k,s}A_{i,k}A_{k,s}^{-1} = A_{i,s}^{-1}A_{i,k}A_{i,s}$.

\item  If $ i<r<k $, $A_{r,s}A_{i,k}A_{r,s}^{-1} =
A_{i,k}^{-1}A_{i,r}^{-1}A_{i,k}A_{i,r}A_{i,k}$.

\item  If $ i<r<k<s$, $A_{r,s}A_{i,k}A_{r,s}^{-1} =
A_{i,s}^{-1}A_{i,r}^{-1}A_{i,s}A_{i,r}
A_{i,k}A_{i,r}^{-1}A_{i,s}^{-1}A_{i,r}A_{i,s}$.
\end{enumerate} Furthermore, these relations are equivalent to the following
relations.
\begin{enumerate}
    \item If either $ s<i$, or $ k<r$, then
    $[A_{i,k}A_{r,s}] = 1$.

 \item  If $ i<k<s$, then
 $[A_{i,k},A_{k,s}^{-1}]= [A_{i,k},A_{i,s}]$.

%% $A_{k,s}A_{i,k}A_{k,s}^{-1} = A_{i,s}^{-1}A_{i,k}A_{i,s}$.
%% $A_{i,k}^{-1}A_{k,s}A_{i,k}A_{k,s}^{-1} = A_{i,k}^{-1}A_{i,s}^{-1}A_{i,k}A_{i,s}$.

\item  If $ i<r<k$, $[A_{r,k}^{-1},A_{i,k}^{-1}]=
[A_{i,k},A_{i,r}]$.

%%$A_{r,k}A_{i,k}A_{r,k}^{-1} =
%%A_{i,k}^{-1}A_{i,r}^{-1}A_{i,k}A_{i,r}A_{i,k}$.

%%$A_{r,k}A_{i,k}A_{r,k}^{-1}A_{i,k}^{-1} =
%%A_{i,k}^{-1}A_{i,r}^{-1}A_{i,k}A_{i,r}$.

\item  If $ i<r<k<s$, $[A_{i,k},A_{r,s}^{-1}] =
[A_{i,k},[A_{i,r},A_{i,s}]]$.

%%$A_{r,s}A_{i,k}A_{r,s}^{-1} =
%%A_{i,s}^{-1}A_{i,r}^{-1}A_{i,s}A_{i,r}
%%A_{i,k}A_{i,r}^{-1}A_{i,s}^{-1}A_{i,r}A_{i,s}$.

%%$A_{i,k}^{-1}A_{r,s}A_{i,k}A_{r,s}^{-1} =
%%A_{i,s}^{-1}A_{i,r}^{-1}A_{i,s}A_{i,r}
%%A_{i,r}^{-1}A_{i,s}^{-1}A_{i,r}A_{i,s}$.

\end{enumerate}
\end{thm}

The well-known values of $A: B_n \to\ \mathrm{Aut}(F_n)$
restricted to $P_n$ given explicitly in \cite{MKS} are listed next
together with some additional related information.

\begin{thm}\label{thm: conjugation action of pure braids} The homomorphism
$A: B_n \to\ \mathrm{Aut}(F_n)$ restricted to $P_n$ is specified
by the following formula:
\[
A(A_{i,k})(x_r)=
\begin{cases}
x_r & \text{if $ r < i$ or $k < r$,}\\
x_i \cdot [x_i,x_k] & \text{if $r = i$,}\\
[x_k, x_i] \cdot x_k & \text{if $r = k$, and}\\
x_r\cdot [x_r,[x_i,x_k]] & \text{if $i< r < k$.}
\end{cases}
\]

In addition, if $ 1 \leq j \leq n$, then
$$\sigma_i^{-1}A_{j,n+1}\sigma_i = A(\sigma_i^{-1})(A_{j,n+1}).$$
The action of $P_n$ given by conjugating the element $A_{j,n+1}$
for $ 1 \leq j \leq n$  by $A_{r,s}$, $ 1 \leq r < s \leq n$ is
given as follows:
$$A_{r,s}^{-1}A_{j,n+1}A_{r,s} = A(A_{r,s}^{-1})(A_{j,n+1})$$

\end{thm}

This theorem specifies the action of the image of $B_n$ in
$\mathrm{Aut}(F_{n+1})$ on the set $A_{1,n+1}, A_{2,n+1}, \cdots,
A_{n,n+1}$. It will be checked below that this action is precisely
that given in the holomorph of a free group. The utility of this
observation is that it then provides easy methods for constructing
certain quotients of braid groups as follows.

Let $\mathbb F$ denote a free group with a choice of generator
given by $x$. Let $G * H$ denote the free product of two groups
$G$, and $H$. There are homomorphisms $$e:\mathrm{Aut}(G) \to\
\mathrm{Aut}(G* \mathbb F),$$ and
$$\chi:G \to\ \mathrm{Aut}(G* \mathbb F)$$
defined as follows.
\begin{defn}\label{defn: E Hol rep}
The homomorphism $e:\mathrm{Aut}(G) \to\ \mathrm{Aut}(G* \mathbb
F)$ is defined on an element $f$ of $\mathrm{Aut}(G)$ by the
formula
\[
[e(f)](z)=
\begin{cases}
f(z) & \text{if $z$ is in $G$, and }\\
z & \text{if $z$ is in $\mathbb F$.}
\end{cases}
\]

The homomorphism $\chi:G \to\ \mathrm{Aut}(G* \mathbb F)$ is
defined on an element $h$ in $G$ by the formula
\[
[\chi(h)](z)=
\begin{cases}
z & \text{if $z$ is in $G$, and }\\
h\cdot z \cdot h^{-1}& \text{if $z$ is in $\mathbb F$.}
\end{cases}
\]

There is an induced function
$$E: \mathrm{Hol}(G) \to\ \mathrm{Aut}(G* \mathbb F)$$
given by the formula $$E(f,h) = e(f) \cdot \chi(h)$$ for $f$ in
$\mathrm{Aut}(G)$, and $h$ in $G$.

\end{defn}

The next classical proposition follows directly with the
verification given below.

\begin{lem}\label{lem:E is a homomorphism}
The following properties hold in $\mathrm{Aut}(G* \mathbb F)$:

\begin{enumerate}
    \item $\chi (h)\cdot e(f)= e(f)\cdot \chi(f^{-1}(h))$,
    \item the function $E: \mathrm{Hol}(G) \to\ \mathrm{Aut}(G* \mathbb F)$ is a
well-defined homomorphism which restricts to $e$ on
$\mathrm{Aut}(G)$, and $\chi$ on $G$ such that $$e(f^{-1})\cdot
\chi(h)\cdot e(f) = \chi(f^{-1}(h)),$$ and
    \item the induced homomorphism $$E: \mathrm{Hol}(G) \to\ \mathrm{Aut}(G* \mathbb
F)$$ is a monomorphism.
\end{enumerate}

\end{lem}

\begin{proof} Notice
that
\[
\chi(h)\cdot e(f)(z)=
\begin{cases}
f(z) & \text{if $z$ is in $G$, and }\\
h\cdot z \cdot h^{-1}& \text{if $z$ is in $\mathbb F$.}
\end{cases}
\]

\[
e(f)\cdot \chi(f^{-1}(h))(z)=
\begin{cases}
f(z) & \text{if $z$ is in $G$, and }\\
h\cdot z \cdot h^{-1}& \text{if $z$ is in $\mathbb F$.}
\end{cases}
\] This relation is exactly that required in the multiplication for
$\mathrm{Hol}(G)$ given by $(f,1)\cdot (1,y)\cdot (f,1)^{-1}= (1,
f^{-1}(y))$, and thus $E$ is a homomorphism which satisfies the
property $$e(f^{-1})\cdot \chi(h)\cdot e(f) = \chi(f^{-1}(h)).$$

To finish, it suffices to check that $E:\mathrm{Hol}(G) \to\
\mathrm{Aut}(G* \mathbb F)$ is a monomorphism. Notice that the
kernel of $E$ projects trivially to $\mathrm{Aut}(G)$, and thus
the kernel is contained in $G$. But $E$ restricted to $G$ is a
monomorphism, and the lemma follows.
\end{proof}

The above constructions now give a direct check that Artin's
representation of the pure braid group is obtained as natural
pull-back. In addition, this identification gives a direct way to
reproduce Artin's result well-known \cite{B} that his
representation is faithful.
\begin{cor}\label{cor: the pullback for the group Pn}
The group $P_{n+1}$ is isomorphic to the group $\zeta_{n+1}$ given
as the pull-back defined by the following cartesian diagram:
\[
\begin{CD}
\zeta_{n+1}  @>{i}>>  P_n  \\
 @VV{I}V          @VV{A|_{P_n}}V \\
\mathrm{Hol}(F_n) @>{j}>>  \mathrm{Aut}(F_n).
\end{CD}
\] Furthermore, Artin's representation of $P_{n+1}$ is given by the
composite
\[
\begin{CD}
P_{n+1}= \zeta_{n+1} @>{I}>>  \mathrm{Hol}(F_n) @>{E}>>
\mathrm{Aut}(F_{n+1}).
\end{CD}
\] Thus Artin's representation is faithful.
\end{cor}

\begin{proof} Let $x_i = A_{i,n+1}$. With this identification, the action of
$P_{n+1}$ on the free group $F_n$ as specified in the extension
given by the holomorph by the formulas $$(f,1)^{-1}\cdot
(1,y)\cdot (f,1)= (1, f^{-1}(y)),$$ and if $ 1 \leq j \leq n$.
That is $$\sigma_i^{-1}A_{j,n+1}\sigma_i =
A(\sigma_i^{-1})(A_{j,n+1})$$ by Theorem \ref{thm: conjugation
action of pure braids}. The next formula follows at once for $ 1
\leq r < s \leq n$:
$$A_{r,s}^{-1}A_{j,n+1}A_{r,s} = A(A_{r,s}^{-1})(A_{j,n+1})= A(A_{r,s}^{-1})(x_j).$$ Thus
these formulas agree with the action given by conjugation in
$P_{n+1}$ by $A_{r,s}$ for $ 1 \leq r < s \leq n$, and so the
pull-back $\zeta_{n+1}$ is isomorphic to $P_{n+1}$.

The next step is to check that Artin's representation $A:P_n \to
\mathrm{Aut}(F_{n})$ is faithful, a statement which is trivially
correct in case $n = 1$. The proof that $A$ is faithful when
restricted to $P_n$ then follows by induction on $n$. That is, the
pull-back of a monomorphism is again a monomorphism while
$P_{n+1}$ is the pull-back of a faithful representation given by
$P_{n+1} \to \mathrm{Hol}(F_n)$ as  in \ref{cor: the pullback for
the group Pn} followed by the injection $E:\mathrm{Hol}(F_n) \to
\mathrm{Aut}(F_{n+1})$ \ref{lem:E is a homomorphism}.

To finish the proof that $A$ is faithful when restricted to $B_n$,
consider the commutative diagram
\[
    \begin{CD}
       B_{n}  @>{A}>> \mathrm{Aut}(F_n) \\
         @VVV          @VV{}V     \\
     \Sigma_{n}  @>{\rho}>> \mathrm{GL}(n, \mathbb Z)
    \end{CD}
\] for which $\rho: \Sigma_{n} \to \mathrm{GL}(n, \mathbb Z)$ is
the natural faithful representation. Thus any element in the
kernel of $B_{n} \to\ \mathrm{Aut}(F_n) $ must be an element of
$P_n$ and is the identity by the previous remarks. The result
follows.
\end{proof}

\section{Quotients of the pure braid group, and the maps $B_n \to\ \mathrm{Aut}(K_n)$ }

Consider the ``reduced free group" due to Milnor\cite{Milnor,
Milnor2} used to analyze ``homotopy string links", and studied in
Habegger-Lin \cite{Habegger-Lin}, \cite{Lin}. The ``reduced free
group" $K_n$ is defined as the quotient of $F_n$ modulo the
relations
$$[x_i,gx_ig^{-1}] = 1$$ where $x_i$ is any generator of $F_n$,
$g$ is any element in $F_n$, and $[x,y]$ denotes the commutator $x
\cdot y \cdot x^{-1} \cdot y^{-1}$. Define $\Lambda_n$ to be the
smallest normal subgroup of $F_n$ containing all of the elements
$[x_i,gx_ig^{-1}]$. Thus $$K_n = F_n / \Lambda_n.$$

This group was rediscovered in a different context \cite{C} to
analyze Barratt's finite exponent conjecture where $K_n$ is the
quotient of the free group on $n$ generators modulo the smallest
normal subgroup containing the simple commutators where at least
one generator appears twice. More precisely, let $\bar \Lambda_n$
denote the smallest normal subgroup of $F_n$ containing all
commutators of the form $[ \cdots [x_{i_1},x_{i_2}]x_{i_3}],\cdots
] x_{i_t}]$ where $ x_{i_j} = x_{i_k}$ for some $ j < k $. The
reduced free group is the quotient
$$K_n = F_n/\bar \Lambda_n.$$ The structure of this group as
well as certain subgroups was analyzed in \cite{C} by considering
the group of units in a non-abelian version of an exterior
algebra.

The next result gives that the action of $B_n$ in
$\mathrm{Aut}(F_n)$ descends to an action on $K_n$, although the
action of the entire automorphism group $\mathrm{Aut}(F_n)$ does
not descend to $K_n$.
\begin{prop}\label{prop: descent}
The action of the group $B_n$ regarded as subgroup of
$\mathrm{Aut}(F_n)$ acting naturally on $F_n$ descends to an
action of $B_n$ on $K_n$ denoted $$\alpha: B_n \to\
\mathrm{Aut}(K_n)$$ and defined by setting
\begin{enumerate}
    \item $\alpha(\sigma_s)(x_j) = x_j$ if $j \neq s, s+1$,
    \item $\alpha(\sigma_s)(x_s) = x_{s+1}$, and
    \item $\alpha(\sigma_s)(x_{s+1}) = x_s[x_{s},x_{s+1}] =
    x_{s+1}^{-1}x_sx_{s+1} $.
\end{enumerate}
\end{prop}

The next technical remark both follows by inspection, and gives an
immediate proof of Proposition \ref{prop: descent}.
\begin{lem}\label{lem:commutator identity}
The following identity is satisfied in any group: $$[x^{-1}yx,
vyv^{-1}] = x^{-1}[y, (xv)y(xv)^{-1}]x $$
\end{lem}

A direct check of this identity is listed next for completeness.
$$[x^{-1}yx, vyv^{-1}]=
x^{-1}\{(y^{-1}(xv)y^{-1}(v^{-1}x^{-1})\}y\{(xv)y(v^{-1})x^{-1}\}x$$
which equals $x^{-1}[y, (xv)y(xv)^{-1}]x$. This lemma is used to
prove Proposition \ref{prop: descent}, a restatement of \ref{prop:
descent of the braid group}.

The proof of Proposition \ref{prop: descent} is given next.
\begin{proof}
It suffices to see that the action $B_n$ on $F_n$ acting naturally
through Artin's representation descends to an action on $K_n$. To
check this point, notice that the relations in $K_n$ are given by
$$[x_s,gx_sg^{-1}]$$ for all $1 \leq s \leq n$ and all $g$ in
$F_n$.

It must be checked that the relations are preserved for the cases
in which
\begin{enumerate}
  \item $x_s$ is replaced by $x_{s+1}$, and
  \item $x_{s+1}$ is replaced by
$x_{s+1}^{-1}x_sx_{s+1} $ with $g$ an arbitrary element of $F_n$.
\end{enumerate} The case $1$ is clear, and the case $2$ is checked
next.

Write $x_{s+1} = x$, and $x_{s} = y$. Then consider
$[x_{s+1},gx_{s+1}g^{-1}] = [x,gxg^{-1}]$. Replace $x$ by
$x^{-1}yx $ to obtain $$[x^{-1}yx,(gx^{-1})y(xg^{-1})].$$ Use the
identity in \ref{lem:commutator identity} $[x^{-1}yx, vyv^{-1}] =
x^{-1}[y,(xv)y(xv)^{-1}]x$ to obtain the relation
$$[x^{-1}yx,(gx^{-1})y(xg^{-1})]= 1.$$ Thus Proposition
\ref{prop: descent} follows.
\end{proof}

With these constructions, there are natural quotients of $P_n$, as
well as an analogous simplicial group obtained from ``homotopy
string links". The reduced braid groups as well as reduced pure
braid groups are defined next.
\begin{defn}\label{defn: reduced braid groups}
Proposition \ref{prop: descent} gives a representation $\alpha:
B_n \to\ \mathrm{Aut}(K_n)$. Define the reduced braid group $rB_n$
as the image of $\alpha$, and reduced pure braid groups $rP_n$ as
the image of $\alpha$ restricted to $P_n$.
\end{defn}

\begin{prop}\label{prop: descent of the braid group}
There is a morphism of group extensions
\[
\begin{CD}
F_n @>{i}>>  P_{n+1} @>{d_i}>>  P_n \\
@VV{}V          @VVV          @VV{\pi}V  \\
K_n @>{i}>> rP_{n+1}   @>{rd_i}>> rP_{n}
\end{CD}
\] where $d_i$ denotes deletion of the $(i+1)$-st strand
with $rd_i$ the induced homomorphism.

In addition, there are morphisms of groups
\[
\begin{CD}
P_{n} @>{s_j}>>  P_{n+1} \\
@VVV          @VV{\pi}V  \\
rP_{n}   @>{rs_j}>> rP_{n+1}
\end{CD}
\] where $s_j$ denotes the degneracy given by``doubling" of the $(j+1)$-st strand
with $rs_j$ the induced homomorphism.

\end{prop}

Thus there are quotients of $P_n$ obtained by replacing $F_q$ by
$K_q$ for all $1 \leq q \leq n-1$ via Definition \ref{defn:
reduced braid groups}. The resulting quotients give groups $rP_n$
of Habegger-Lin, and Milnor which arises from the notion of link
homotopy as described crudely in the next paragraph.

Consider the space of ordered $n$-tuples of  smooth maps of $S^1$
in $\mathbb R^3$ having disjoint images. Milnor defines an
equivalence relation of "link homotopy" by allowing strands from
the same component of a link to pass over each other
\cite{Milnor,Milnor2}. Lin describes the infinitesimal
link-homotopy relations on page $5$ of \cite{Lin}.

For example, consider the trivial $n$-component link
$\mathcal{L}_n$ in $\mathbb R^3$. Thus the fundamental group of
the complement $\mathbb R^3 - \mathcal{L}_n$ is isomorphic to
$F_n$ by some choice of isomorphism. Let $\beta$ denote a simple
closed curve in the complement which represents an element $b$ in
$F_n$. Then both Milnor, and Habegger-Lin prove that the link
$\mathcal{L}_n \cup \beta$ is link homotopically trivial if and
only if the element $b$ projects to the identity in $K_n$.
\cite{Milnor,Milnor2, Habegger-Lin, Lin}.

\begin{thm}\label{thm: Free to K }
The natural quotient map $\rho: F_n \to\ K_n$ prolongs to a
morphism of simplicial groups $$\rho: F[S^1] \to\ K[S^1]$$ for
which $K[S^1]$ in degree $n$ is $K_{n}$. The kernel of $\rho$,
$\Gamma[S^1]$, is a free group in each degree.

\begin{enumerate}
  \item The map $\rho$ is a
surjection of simplicial groups, and is thus a fibration with
fibre given by a simplicial group which is free in each degree.
  \item The $n$-th homotopy group of $K[S^1]$, $\pi_nK[S^1]$,
  is isomorphic to $\oplus_{(n-1)!}\mathbb Z$, the $(n-1)$-st
  homology group of $P_n$ given by $\mathrm{Lie}(n)$.
  \item The induced map $\pi_*(\rho): \pi_*(F[S^1]) \to\ \pi_*(K[S^1])$
  is an isomorphism in dimension $1$, and $2$. This map is zero in
  degrees greater than $2$.
  \item If $n>2$, here is a short exact sequence of homotopy groups
$$0 \to\ \mathrm{Lie}(n+1) \to\ \pi_{n} \Gamma[S^1] \to\ \pi_n(F[S^1]) \to\ 0. $$
\end{enumerate}
\end{thm}

\begin{proof}
The quotient maps $\rho: F_n \to\ K_n$ commute with the face and
degeneracies by \ref{prop: descent of the braid group}. That
$\rho$ is a surjection of simplicial groups is a consequence of
\ref{prop: descent of the braid group}. In addition, the homotopy
groups of $K[S^1]$ were determined in \cite{C}. The theorem
follows at once from standard properties of simplicial groups
applied to $K[S^1]$ \cite{C,W}.
\end{proof}

The above remarks provide a comparison of the group theory
associated to link homotopy to the group theory associated to the
homotopy groups of the $2$-sphere. In addition, there is a functor
from simplicial groups to pro-groups which when specialized to
$K[S^1]$ gives a group $H_{\infty}$ \cite{C2} which is filtered
with associated graded given by $$E_0^*(H_{\infty}) = \oplus_{n
\geq 1} Lie(n).$$ The framework above provides other natural
quotients of $F_n$, as well as $P_n$ which will be addressed
elsewhere.

\section{ On link homotopy and related homotopy groups }
\label{sec: On link homotopy and related homotopy groups   }

The purpose of this section is to describe some connections of
free groups, braid groups, homotopy groups of the $2$-sphere. The
subsequent sections will illustrate one connection between "link
homotopy", Artin's pure braid groups, and the homotopy groups of
certain choices of simplicial groups.

There is a braid invariant obtained from homotopy groups as
follows. Consider the embedding $\Theta_n:F_n \to\ P_{n+1}$. Let
$C_n$ denote the chains in degree $n$ for $F[S^1]$. That is the
normal subgroup given by the intersection of the kernels of $d_1,
d_2, \cdots, d_n$ (excluding $d_0$).Thus $P_{n+1}$ is equal to a
disjoint union of left cosets:
$$P_{n+1} = \amalg_{\alpha \epsilon S} x_{\alpha}C_n$$ where
$$\{x_{\alpha}|\alpha \epsilon S \}$$ is a complete set of
distinct left coset representatives for $C_n$ in $P_{n+1}$. Thus
given a pure braid $\gamma$, assign the homotopy element given by
the homotopy class of $x$ where $\gamma = \alpha \cdot x$. As an
example, this element is analyzed for three-stranded braids.

\begin{exm} Let $E_0^n(K_q)$  denote the $n$-stage of the
descending central series for $K_q$ modulo the $n+1$-st stage.
Then $E_0^n(K_3)$ has a basis as follows \cite{C2}.
\begin{enumerate}
\item If $n=1 $, a basis is $x_1$, $x_2$, and $x_3$. \item If
$n=2$, a basis is $[x_1,x_2]$, $[x_1,x_3]$, and $[x_2,x_3]$. \item
If $n = 3$, a basis is $[[x_1,x_2],x_3]$, and $[[x_1,x_3],x_2]$.
\end{enumerate}
A basis for $E_0^t(K_n)$ is
$[[x_{i_1},x_{i_\tau(2)}],x_{i_\tau(3)}], \cdots ],x_{i_\tau(t)}]$
for all sequences $ 1\leq i_1 < i_2 < i_3< \cdots < i_t \leq n $
where $\tau$ runs over all elements in the symmetric group
$\Sigma_{t-1}$.
\end{exm}

A non-bounded cycle in $F_2$ is given by $[x_1,x_2]$. This cycle
represents a choice of generator of $\pi_2 \Omega S^2 = \mathbb Z$
given by the classical Hopf map. In addition, the braid closure of
$[x_1,x_2]$ gives the Borromean rings. The computations are
direct, and left as an exercise.

\section{On looping $\mathrm{AP}_*$}

The next theorem was proven in \cite{CW} while a detailed proof is
included here for completeness.
%%Theorem
\begin{thm} \label{theorem:Looping braids}
The loop space ( as simplicial groups ) of the simplicial group
$\mathrm{AP}_*$, $\Omega(\mathrm{AP}_*)$, is isomorphic to
$F[\Delta[1]]$ as a simplicial group. Thus $\mathrm{AP}_*$ is
contractible, and the realization of the simplicial set
$\mathrm{AP}_*/ F[S^1]$ is homotopy equivalent to $S^2$.
\end{thm}

\begin{remark}
The long exact homotopy sequence of the fundamental fibration
sequence due to Fadell, and Neuwirth \cite{FN} can be regarded as
determining the homotopy type of the loop space for the simplicial
group obtained from the pure braid groups.
\end{remark}

\begin{proof}
John Moore gave a definition for the loop space $\Omega \Gamma_*$
of a reduced simplicial group $\Gamma_*$ ( a simplicial group for
which  $\Gamma_0$ consists of a single element) \cite{Moore}. This
procedure corresponds to the topological notion of looping in the
sense that the loop space of the geometric realization of
$\Gamma_*$ is homotopy equivalent to the geometric realization of
$\Omega \Gamma_*$. This process of looping a simplicial group is
described next.

Define a simplicial group $E\Gamma_*$ where the group $E\Gamma_n$
in degree $n$ is given by the group $\Gamma_{n+1}$ with face, and
degeneracies given by "shifting down by one". Then $\Omega
\Gamma_*$, the looping of $\Gamma_*$ in degree $n$, is defined to
be the kernel of the map $$d_0: E\Gamma_n \to\ \Gamma_n$$ thus
$$\Omega \Gamma_n =  \mathrm{ker}[d_0: \Gamma_{n+1} \to\ \Gamma_n ].$$

In the special case that $\Gamma_n$ is $ P_{n+1}$, then
$E\Gamma_n$ is $ P_{n+2}$. Recall that $P_{n+1}$ is generated by
symbols $A_{i,j}$ for $1 \leq i < j \leq n+1$. Furthermore, the
map $d_0: E\Gamma_n \to\ \Gamma_n$ is induced by the projection
map $p_*: \pi_1(F(\mathbb R^2,n+2)) \to\ \pi_1(F(\mathbb
R^2,n+1))$ where $p: F(\mathbb R^2,n+2) \to\ F(\mathbb R^2,n+1)$
is the map given by projection to the first $n+1$ coordinates.

The fibre of the map $p$ is $\mathbb R^2- Q_{n+1}$, and thus the
kernel of $p_*$ is isomorphic to $F_{n+1}$. It will be shown below
that a choice of generators for this kernel, regarded as a
simplicial set is $\Delta[1]$. Together with Theorem $1.2$ giving
that $\Theta: F[S^1] \to\ \mathrm{AP}_*$ is an embedding, the
result follows.

The kernel of the map $d_0: P_{n+2} \to\ P_{n+1}$ which is given
by deleting the last coordinate is given by the free group with
generators $P_{n+2,j}$ for $1 \leq j < n+2$. A basis for the
kernel is given by $\{P_{n+2,j}| 1 \leq j \leq n+1 \}$.

Furthermore, the simplicial $1$-simplex is given by $\langle
0^i,1^{n+1-i}\rangle$ in degree $n$ with $ 0 \leq i \leq n+1 $
with the single relation that $\langle 0^{n+1}\rangle= 1$ in the
simplicial group $F[\Delta[1]$. There is a dimension-wise morphism
of groups
$$\Psi: F[\Delta[1]] \to\ \Omega(\mathrm{AP}_*)$$ defined by
sending $\langle 0^i,1^{n+1-i}\rangle$ to $P_{n+2,i}$ for $ 0 < i
< n+2 $. A direct check gives that this dimension-wise
homomorphism is compatible with the face, and  degeneracies. Since
each map is an isomorphism, the map $\Psi$ is an isomorphism of
simplicial groups. The theorem follows.
\end{proof}

\section{On embeddings of residually nilpotent groups}

Let $\rho: \pi \to\ G$ be a homomorphism between discrete groups.
Let $\Gamma^n(\pi) = \Gamma^n$ denote the $n$-th stage of the
descending central series for $\pi$. That is
\begin{enumerate}
  \item  $\Gamma^1(\pi) = \pi$, and inductively
  \item  $\Gamma^{n+1}(\pi) = [\Gamma^n,\pi]$,
  the subgroup generated by
  commutators $[\cdots [h_1,h_2]h_3]\cdots]h_q]$ with
  $h_i$ in $\pi$, and $ q \geq n+1$.
\end{enumerate} Define the associated graded
$$ E_0^n(\pi)= \Gamma^{n}(\pi)/\Gamma^{n+1}(\pi),$$ and
$$E_0^*(\pi) = \oplus_{ n \geq 1} E_0^n(\pi).$$ It is
a classical and easily checked fact that the commutator
$$[-,-]:\pi \times \pi \to\ \pi$$ given by $$[x,y] =
x^{-1}\cdot y^{-1}\cdot x\cdot y$$ induces the structure of Lie
algebra on $E_0^*(\pi)$. A group homomorphism  $\rho:\pi \to\ G$
preserves the stages of the descending central series. Thus there
is an induced morphism of associated graded Lie algebras
$$E_0^*(\rho): E_0^*(\pi) \to\ E_0^*(G).$$

\begin{defn}\label{defn: residually nilpotent}
A discrete group $\Gamma$ is said to be residually nilpotent group
if $$\bigcap_{i \geq 1}\Gamma^i(\pi) = \{\mathrm{identity}\}$$
where $\Gamma^i(\pi)$ denotes the $i$-th stage of the descending
central series for $\pi$.
\end{defn}

\begin{thm}\label{thm: one to one}
Assume that $\pi$ is a residually nilpotent group. Let $$\rho: \pi
\to\ G$$ be a homomorphism of discrete groups such that the
morphism of associated graded Lie algebras
$$E_0^*(\rho): E_0^*(\pi) \to\
E_0^*(G)$$ is a monomorphism. Then $\rho$ is a monomorphism.
\end{thm}

\begin{proof}
Let $x$ denote a non-identity  element in the kernel of $\rho$.
Since $\pi$ is residually nilpotent, there exists a natural number
$n$ such that the element $x$ is in $\Gamma^n( \pi )$ and not in
$\Gamma^{n+1}( \pi )$. But then $x$ projects to a non-identity
element in $E_0^n(\pi)$, and thus has non-trivial image in
$E_0^n(G)$ contradicting the fact that $x$ is a non-identity
element in the kernel of $\rho$.
\end{proof}

If $F[S]$ is a free group generated by the set $S$, then $F[S]$ is
residually nilpotent \cite{MKS}, and the next corollary follows at
once.
\begin{cor} \label{prop:embeddings}
If $F[S]$ is a free group generated by the set $S$, and
$$E_0^*(\rho): E_0^*(F[S]) \to\ E_0^*(G)$$ is a
monomorphism of Lie algebras, then $\rho: F[S] \to\ G$ is a
monomorphism of groups. If $\rho$ is assumed to be an epimorphism,
then it is an isomorphism.
\end{cor}

\begin{proof}
If $f$ induces an monomorphism on $E_0^*(F[S])$, then $f$ is a
monomorphism by the previous theorem. It suffices to note that $f$
induces a surjection.  The lemma follows.
\end{proof}

This approach for testing whether a homomorphism $\rho$ is
faithful is suited for residually nilpotent groups $\pi$ such as
free groups, and pure braid groups $P_{n+1}$. A single case is
addressed here in section $8$ here. One, possibly useful case, is
that of the Gassner representation.

\section{The simplicial structure for $\mathrm{AP}_*$}
\label{sec:The simplicial structure for $AP_*$}

Recall that $P_{n+1}$ is generated by symbols $A_{i,j}$ for $1
\leq i < j \leq n+1$. Artin's relations are listed in \cite{MKS}
as well as in section $3$ here.

The face operations in the simplicial group $\mathrm{AP}_*$ are
defined as follows:
%%%%%%%    conventions as in the Skye paper      %%%%%
\[
d_t(A_{i,j})  =
\begin{cases}
A_{i-1,j-1}
  & \quad   \text{if $ t+1 < i $, \ } \\
1
  & \quad   \text{if t +1 = i, \ } \\
A_{i,j-1}
  & \quad   \text{if $i < t + 1 < j$ , \ } \\
1
  & \quad   \text {if t + 1 = j, \ } \\
A_{i,j}
  & \quad   \text{if $t + 1 > i$ . \ } \\
\end{cases}
\]

The degeneracy operations are defined as follows:
\[
s_t(A_{i,j}) =
\begin{cases}
A_{i+1,j+1}
  & \quad   \text{if $t+1 < i$ , \ } \\
A_{i,j+1}\cdot A_{i+1,j+1}
  & \quad   \text{if t+1 = i , \ } \\
A_{i,j+1}
  & \quad   \text{if $i < t + 1 < j$ , \ } \\
A_{i,j}\cdot A_{i,j+1}
  & \quad   \text{if t + 1 = j , \ } \\
A_{i,j}
  & \quad   \text{if $t + 1 > j$. \ } \\
\end{cases}
\]

That the face and degeneracies preserve the defining relations is
verified directly. The details are omitted. These operations give
a convenient method for describing the behavior of $\Theta_n$ in
the next section.

The face, and degeneracy operations for $F[S^1]$ are prescribed as
follows.
\[
d_i(x_j) =
\begin{cases}
x_j
  & \quad   \text{if $j < i$ , \ } \\
1
  & \quad   \text{if i = j , and \ } \\
x_{j-1}
  & \quad   \text{if $j > i$ . \ } \\
\end{cases}
\]

\section{The Lie algebra associated to the descending central series for $P_{n+1}$}
\label{sec:The Lie algebra associated to the descending central
series for $P_{n+1}$}

The next theorem gives the structure of the Lie algebra arising
from the descending central series for $P_{k}$ which was analyzed
in work of T. Kohno \cite{K,K1}, Falk, and Randell \cite{FR}, and
Xicot\'encatl \cite{Xico}. Let $B_{i,j}$ denote the projections of
the $ A_{i,j}$ to $E_0^*(P_k)$.

\begin{thm} \label{thm: Lie algebras}
The Lie algebra obtained from the descending central series for
$P_{k}$ is given by $\mathcal{L}_k $  the free Lie algebra
generated by elements $B_{i,j}$ with $1 \leq i < j \leq k$, modulo
the infinitesimal braid relations:
\begin{description}
\item[(i)] $[B_{i,j}, B_{s,t}] =0 $ if $\{ i,j\} \cap \{ s,t\} =
\phi,$ \item[(ii)] $[ B_{i,j}, B_{i,t} + B_{t,j}] = 0$ if $1 \leq
i < t < j \leq k $, and \item[(iii)] $[B_{t,j}, B_{i,j} + B_{i,t}]
= 0$ if $1 \leq i < t < j \leq k$.
\end{description}

Furthermore there is a split short exact sequence of Lie algebras
\[
\begin{CD}
0 \to\ E_0^*(F_n) @>{ E_0^*(i) }>> E_0^* (P_{n+1})
@>{E_0^*(d_n)}>> E_0^*(P_n)\to\ 0
\end{CD}
\] where $E_0^*(F_n)$ is the free Lie algebra generated by
$B_{i,n+1}$ for $1 \leq i < n+1$. In addition, $E_0^* (P_{n+1})$
is additively isomorphic to $ E_0^*(P_n) \oplus E_0^*(F_n).$
\end{thm}

A related direct computation which is used below follows.
\begin{prop} \label{prop: infinitesimal braid relations again}
\begin{enumerate}
    \item If $ 1 \leq j < t < k \leq  n+1 $, then
    $$[B_{j,k} + B_{t,k},B_{j,t}]=0,$$ and
    $$[B_{j,t},B_{t,k}] = [B_{t,k},B_{j,k}].$$
    \item If $ 1 \leq r < m < n+1 $, then
    $$[\sum_{1 \leq i \leq n} B_{i,n+1},
    \sum_{ 1 \leq j \leq r } B_{j,m}] = 0.$$

\end{enumerate}
\end{prop}

\begin{proof}
Assume that $ 1 \leq j < t < k \leq  n+1 $, and consider the
infinitesimal braid relations as follows:
\begin{itemize}
  \item $[B_{j,t} + B_{j,k},B_{t,k}]= 0$
  \item $[B_{j,t} + B_{t,k},B_{j,k}]= 0$
\end{itemize} Thus $[B_{j,k} + B_{t,k},B_{j,t}]= [B_{t,k}, B_{j,k}]+
[B_{j,k},B_{t,k}] = 0$.

In addition, $[B_{j,t}+B_{j,k},B_{t,k}] = 0$, and thus
$[B_{j,t},B_{t,k}] = [B_{t,k},B_{j,k}]$.

Consider $[\sum_{1 \leq i \leq n} B_{i,n+1}, B_{j,m}]$ for $m <
n+1 $. Since $[B_{i,n+1},B_{j,m}] = 0 $ if $\{i,n+1\} \cap \{j,m\}
= \phi$, it follows that $[\sum_{1 \leq i \leq n} B_{i,n+1},
B_{j,m}] = [B_{m,n+1} + B_{j,n+1},B_{j,m}]= 0$ by the
infinitesimal braid relations.

The proposition follows.
\end{proof}

\section{On $\Theta_n: F_{n} \to\ P_{n+1}$}\label{sec:On Theta n}

The proof that the map $\Theta_n: F_{n} \to\ P_{n+1}$ is a
monomorphism depends on the structure of certain Lie algebras
given in this section. These Lie algebras arise from passage to
the associated graded for the descending central series filtration
of a discrete group arising in the commutative diagram of groups
\[
\begin{CD}
F_{n}  @>{\Theta_n}>>  P_{n+1}  \\
@VV{d_n}V             @VV{d_n}V       \\
F_{n-1}    @>{\Theta_{n-1}}>>  P_{n}.
\end{CD}
\]

This diagram together with induction on $n$ is used to prove the
following theorem.
\begin{thm}\label{thm: embeddings of Lie algebras}
The maps $\Theta_n :F[x_1,x_2, \cdots, x_n] \to\ P_{n+1}$ on the
level of associated graded Lie algebras
$$E_0^*(\Theta_n) :E_0^*(F[x_1,x_2, \cdots, x_n]) \to\ E_0^*(P_{n+1})$$
are monomorphisms. Thus, by Theorem \ref{thm: one to one}, the
maps $\Theta_n$ are monomorphisms.
\end{thm}

The hypothesis that $\Theta_1$ is an isomorphism gives the initial
step in an induction with the assumption that
$E_0^*(\Theta_{n-1})$ is an embedding. To carry out the inductive
step, notice that there is a commutative diagram of morphism of
Lie algebras
\[
\begin{CD}
E_0^*(F_{n})  @>{E_0^*(\Theta_n)}>>  E_0^*(P_{n+1})  \\
@VV{E_0^*(d_n)}V             @VV{E_0^*(d_n)}V \\
E_0^*(F_{n-1})    @>{E_0^*(\Theta_{n-1})}>>  E_0^*(P_{n}).
\end{CD}
\]

Most of this section gives explicit results concerning these Lie
algebras which are used to prove Theorem \ref{thm: embeddings of
Lie algebras} in the next section. Thus recall P. Hall's classical
result that
$$E_0^*( F_{n}) = E_0^*(F[x_1,x_2, \cdots,
x_n])$$ is isomorphic to the free Lie algebra over the integers
$\mathbb Z$
$$L[x_1,x_2, \cdots, x_n]$$ where the $x_i$'s in the free Lie
algebra are the projections of the analogous elements in the group
$F_n$. Furthermore, the kernel of the projection map
$$\pi: L[x_1,x_2, \cdots, x_n] \to\ L[x_2, \cdots, x_n]$$
defined by
\[
\pi(x_i)=
\begin{cases}
0 & \text{if $i=1$,}\\
x_i & \text{if $i > 1$.}
\end{cases}
\]
is given by $L[S_n]$ where $S_n$ is the set
$$\{x_1, \mathrm{ad}(x_{i_1})\mathrm{ad}(x_{i_2})\cdots \mathrm{ad}(x_{i_p})(x_1)\}$$
for all $p \geq 1$ with $i_j > 1$. Thus the Lie algebra kernel of

$$E_0^*(d_n): E_0^*(F_{n}) \to E_0^*(F_{n-1})$$
is given by $L[S_n]$.

It is not the case that the exact sequence of groups $ 1 \to\
\mathrm{ker}(d_n) \to\  F_{n}  \to\ F_{n-1} \to\ 1$ induces a
short exact sequence of Lie algebras after passage to the
sub-quotients of the descending central series. The details of
this last assertion are not included.

It is convenient to give define additional elements specified in
the following formulae.
\begin{enumerate}
  \item $\Lambda_n = B_{1,n+1} + B_{2,n+1}+ \cdots + B_{n,n+1}$, and
  \item $\gamma_q(n) = -\Sigma_{n-q+2 \leq i \leq n} B_{i, n+1}$ if
  $2 \leq q \leq n $.
\end{enumerate}
Some values of $\Theta(x_i)$ are required next, and are recorded
in the next theorem which is a direct computation.

\begin{thm}\label{thm: values of theta}
The map $\Theta_n :F[x_1,x_2, \cdots, x_n] \to\ P_{n+1}$ satisfies
the following formula:

$$E_0^*(\Theta_n)(x_q) = \Sigma_{1 \leq i \leq n-q+1 < j \leq n+1}
B_{i,j} $$ on the level of associated graded Lie algebras
$E_0^*(P_{n+1})$ for $ 1 \leq q \leq n$. In addition, the
following formulas hold.
\begin{enumerate}

\item $E_0^*(\Theta_n)(x_1) = \Lambda_n$.

\item If $j < n+1$, $[\Lambda_n, B_{i,j}] = 0$.

\item $E_0^*(\Theta_n)(x_q) = \Lambda_n + \gamma_q(n) + \Sigma_{1
\leq i \leq n-q+1<  j \leq n} B_{i,j}$.

\item If $2 \leq p \leq q$, then $[\gamma_p(n), \Theta(x_q(n))]=
0$.

\end{enumerate}
\end{thm}

\begin{proof}
The formula
$$E_0^*(\Theta_n)(x_q) = \Sigma_{1 \leq i \leq n-q+1,n-q+2
\leq j \leq n}B_{i,j}$$ is that given by the degeneracies as
described in section $6$, the simplicial structure for $A P_*$.
The next $4$ formulae are proven next with formula $1$ given by
the definition of $\Lambda_n$.

With the condition $s < t < n+1$, notice that $[B_{i,n+1},B_{s,t}]
= 0$ for $i \neq s,t$. Thus
$$[\Lambda_n,B_{s,t}] = [\Sigma_{1 \leq i \leq
n}B_{i,n+1},B_{s,t}],$$ and
$$[\Lambda_n,B_{s,t}] = [B_{s,n+1}+B_{t,n+1},B_{s,t}]$$ which is $0$
by the infinitesimal braid relations in Proposition \ref{prop:
infinitesimal braid relations again}. Formula $2$ follows.

Since $$E_0^*(\Theta_n)(x_q) = \Sigma_{1 \leq i \leq n-q+1,n-q+2
\leq j \leq n}B_{i,j},$$ and $$\Lambda_n + \gamma_q(n) = \Sigma_{1
\leq i \leq n-q+1}B_{i,n+1},$$ formula $3$ that
$E_0^*(\Theta_n)(x_q) = \Lambda_n + \gamma_q(n) + \Sigma_{1 \leq i
\leq n-q+1 < j \leq n} B_{i,j}$ follows directly.

To work out part $4$, consider $$[\gamma_p(n), \Theta(x_q)]=
-[\Sigma_{n-p+ 1 \leq i \leq n} B_{i, n+1}, \Sigma_{1 \leq i \leq
n-q+1 < \leq j \leq n+1} B_{i,j}].$$ Then expand each term for
$n-p+ 1 \leq i \leq n$ given by
$$[B_{i, n+1},
\Sigma_{1 \leq i \leq n-q+1 < \leq j \leq n+1} B_{i,j}].$$ There
are two cases to check.
\begin{enumerate}
  \item If $i < n-q+1$, then $[B_{i, n+1},
\Sigma_{1 \leq i \leq n-q+1 < \leq j \leq n+1} B_{i,j}]= 0$ by the
infinitesimal braid relations.
  \item If $i = n-q+1$, then $[B_{n-q+1, n+1},
\Sigma_{1 \leq i \leq n-q+1 < \leq j \leq n+1} B_{i,j}]= 0$ by the
infinitesimal braid relations.
\end{enumerate}

The theorem follows.
\end{proof}

Some additional properties concerning the degeneracies are given
next follow at once.
\begin{thm}\label{thm: values of degeneracies and theta}
The map $\Theta_n :F[x_1,x_2, \cdots, x_n] \to\ P_{n+1}$ satisfies
the following formula:
\begin{enumerate}

\item For any fixed $0 \leq j \leq n$, $s_j\Lambda_{n} =
\Lambda_{n+1}$.

\item
\[
s_j(\gamma_q(n)) =
\begin{cases}
\gamma_q(n+1) & \text{if $j < n+1-q$,}\\
\gamma_{q+1}(n+1) & \text{if $j \geq n+1 - q$.}
\end{cases}
\]

\item $[\gamma_3(3),\Theta(x_2)]= [\gamma_3(3),\Lambda_3]+
[\Lambda_3,\gamma_2(3)] + 2[\gamma_3(3),\gamma_2(3)]$.

\item For any fixed $2 \leq i< j \leq n$, there are sequences of
degeneracies $s(I,J)$ such that
\begin{itemize}
  \item $s(I,J)(x_2) = x_i$,
  \item $s(I,J)(x_3) = x_j$,
  \item $s(I,J)([\gamma_3(3),E_0^*(\Theta_n)(x_2)] =
  [\gamma_j(n),E_0^*(\Theta_n)(x_i)]$.
\end{itemize}

\item For any fixed $2 \leq i< j \leq n$,
$$[\gamma_j(n),E_0^*(\Theta_n)(x_i)] =
  [\gamma_j(n),\Lambda_n]+ [\Lambda_n,\gamma_i(n)] +
 2[\gamma_j(n),\gamma_i(n)].$$
%%6
\item If $2 \leq p \leq q $, then
$$[\gamma_p(n), \Theta_n(x_q(n))]= 0.$$

\item If $2 \leq q < p$, then
$$[\gamma_p(n), \Theta_n(x_q(n)) - 2\gamma_q(n)]=
[\gamma_p(n),\Lambda_n]+ [\Lambda_n,\gamma_q(n)].$$
\end{enumerate}
\end{thm}

\begin{proof}
The degeneracies are described in section \ref{sec:The simplicial
structure for $AP_*$} which gives the simplicial structure for
$\mathrm{AP}_*$. Recall that $\Lambda_n = \Sigma_{1 \leq i \leq
n}B_{i,n+1}$. That $s_j\Lambda_{n} = \Lambda_{n+1}$, and formula
$1$ follows at once. Recall that $\gamma_q(n) = B_{n+2-q,n+1}+
B_{n+3-q,n+1}+ \cdots + B_{n,n+1}$. Thus if $j \leq n+1 -q$, then
$s_j\gamma_q(n) = B_{n+3-q,n+2}+ B_{n+4-q,n+2}+ \cdots +
B_{n+1,n+2}$. Furthermore,  if $n+2 -q \leq j$, then
$s_j\gamma_q(n) = B_{n+2-q,n+2}+ B_{n+3-q,n+2}+ \cdots +
B_{n+1,n+2}$, and formula $2$ follows.

To check formula $3$, notice that the following hold.
\begin{enumerate}
  \item By definition, $\gamma_3(3) = -(B_{2,4}+ B_{3,4})$.
  \item By definition, $E_0^*(\Theta_3)(x_2) = B_{1,3}+ B_{2,3}+ B_{1,4}+
  B_{2,4}$.
  \item Thus $[\gamma_3(3),E_0^*(\Theta_3)(x_2)] = -[B_{2,4}+ B_{3,4},B_{1,3}+ B_{2,3}+
  B_{1,4}+ B_{2,4}]$ which equals
  \item $-[B_{2,4},B_{1,3}+ B_{2,3}+
  B_{1,4}+ B_{2,4}]$ by the infinitesimal braid relations.
  \item Thus $[\gamma_3(3),E_0^*(\Theta_3)(x_2)] = -[B_{2,4}, B_{2,3}+
  B_{1,4}]$ by the infinitesimal braid relations.
  \item Furthermore, $[B_{2,4}, B_{2,3}+ B_{3,4}] = 0$ by the
  infinitesimal braid relations, and thus
  \item $[\gamma_3(3),\Theta(x_2)] = [B_{2,4}, B_{3,4}] -
  [B_{2,4}, B_{1,4}] = [B_{2,4}, B_{3,4}] - [B_{2,4}, \Lambda_3 -
  B_{3,4}]$.
  \item Substituting $\gamma_2(3)= -B_{3,4}$, $\gamma_3(3)= -
  (B_{2,4}+ B_{3,4}$, and $\gamma_2(3) - \gamma_3(3) = B_{2,4}$
  gives $[\gamma_3(3),\Theta(x_2)] = [\gamma_3(3),\Lambda_3]+
[\Lambda_3,\gamma_2(3)] + 2[\gamma_3(3),\gamma_2(3)]$. Thus
formula $3$ follows.
\end{enumerate}

Formula $4$ clearly follows from the degeneracies.

To check formula $5$, notice that
$$[\gamma_j(n),E_0^*(\Theta_n)(x_i)] =
s(I,J)([\gamma_3(3),E_0^*(\Theta_n)(x_2)],$$
$$s(I,J)([\gamma_3(3),E_0^*(\Theta_3)(x_2)] =
S(I,J)([\gamma_3(3),\Lambda_3]+ [\Lambda_3,\gamma_2(3)] +
2[\gamma_3(3),\gamma_2(3)]),$$ and
$$S(I,J)([\gamma_3(3),\Lambda_3]+[\Lambda_3,\gamma_2(3)]+2[\gamma_3(3),
\gamma_2(3)])$$ equal to
$$
[\gamma_j(n),\Lambda_n]+[\Lambda_n,\gamma_i(n)]+2[\gamma_j(n),\gamma_i(n)].$$

Formula $6$ is formula $4$ in Theorem \ref{thm: values of theta}.

By formula $5$, For any fixed $2 \leq i< j \leq n$,
$$[\gamma_j(n),E_0^*(\Theta_n)(x_i)] =
  [\gamma_j(n),\Lambda_n]+ [\Lambda_n,\gamma_i(n)] +
 2[\gamma_j(n),\gamma_i(n)].$$ Thus
$$[\gamma_j(n),E_0^*(\Theta_n)(x_i) -2\gamma_i(n) ] =
  [\gamma_j(n),\Lambda_n]+ [\Lambda_n,\gamma_i(n)],$$
  and formula $7$ follows.

\end{proof}

\section{On the proof of Theorem \ref{thm: embeddings of Lie algebras}}\label{sec: lie algebras}

Recall the morphism of Lie algebras
\[
\begin{CD}
E_0^*(F_{n})  @>{E_0^*(\Theta_n)}>>  E_0^*(P_{n+1})  \\
@VV{E_0^*(d_n)}V             @VV{E_0^*(d_n)}V \\
E_0^*(F_{n-1})    @>{E_0^*(\Theta_{n-1})}>>  E_0^*(P_{n}).
\end{CD}
\] Taking Lie algebra kernels results in a commutative diagram

\[
\begin{CD}
L[x_1^{E_0^*(F_{n-1})}]  @>{E_0^*(\Theta_n)}>>
L[B_{1,n+1},B_{2,n+1},\cdots,B_{n,n+1}] \\
@VV{}V             @VV{}V \\
E_0^*(F_{n})  @>{E_0^*(\Theta_n)}>>  E_0^*(P_{n+1})  \\
@VV{E_0^*(d_n)}V             @VV{E_0^*(d_n)}V \\
E_0^*(F_{n-1})    @>{E_0^*(\Theta_{n-1})}>>  E_0^*(P_{n})
\end{CD}
\] for which
\begin{itemize}
  \item the Lie algebra kernel of
  $E_0^*(d_n):E_0^*(F_{n}) \to\ E_0^*(F_{n-1})$ is
  $L[x_1^{E_0^*(F_{n-1})}]$, and
  \item the Lie algebra kernel of
  $E_0^*(d_n):E_0^*(P_{n+1}) \to\ E_0^*(P_{n})$ is
  $$L[B_{1,n+1},B_{2,n+1},\cdots,B_{n,n+1}].$$
\end{itemize}

The inductive hypothesis is that $E_0^*(\Theta_{n-1})$ is a
monomorphism. To finish, it suffices to check that the induced map
of Lie algebras
\[
\begin{CD}
L[x_1^{E_0^*(F_{n-1})}]  @>{E_0^*(\Theta_n)}>>L[B_{1,n+1},
B_{2,n+1},\cdots,B_{n,n+1}]
\end{CD}
\] is a monomorphism. This will be checked using the computations
of section \ref{sec:The Lie algebra associated to the descending
central series for $P_{n+1}$}.

The next step is to consider the morphism of Lie algebras
$$p: L[B_{1,n+1}, B_{2,n+1},\cdots,B_{n,n+1}] \to\
L[B_{2,n+1},\cdots,B_{n,n+1}]$$ defined by the formula
\[
p(B_{j,n+1})=
\begin{cases}
B_{j,n+1} & \text{if $j > 1$,}\\
-(B_{2,n+1}+ B_{3,n+1} + \cdots + B_{n,n+1}) & \text{if $j=1$.}
\end{cases}
\] Notice that $p( \Sigma_{1 \leq j \leq n}B_{j,n+1}) =0,$ and so
$$p(\Lambda_n) = 0.$$ A remark is appropriate here. With the
exception of the map $p$, all of the above maps are morphisms of
simplicial Lie algebras. The map $p$ fails to preserve one face
operation, and is thus not a morphism of simplicial Lie algebras.

Let $B$ equal  $E_0^*(F_{n-1})$. The image of
$$E_0^*(\Theta_n): L[x_1^{B}]\to\ L[B_{1,n+1}, B_{2,n+1},
\cdots,B_{n,n+1}]$$ is in the Lie ideal generated by $\Theta(x_1)
= \Lambda_n$, and thus is in the Lie algebra kernel of $p$. Hence
the map $\Theta$ restricts to a map of Lie algebras $$
E_0^*(\Theta_n): L[x_1^{B}]\to\ L[\Lambda_n^{C}]$$ where $C=
L[B_{2,n+1},B_{3,n+1},\cdots,B_{n,n+1}]$. In the work which is
given below, the notation $x_j = x_j(n)$, and $\gamma_j =
\gamma_j(n)$ is used when the integer $n$ is clear from the
context.

The structure of the Lie algebras  $L[x_1^{B}]$, and
$L[\Lambda_n^{C}]$ are used below. Recall that the abelianization
of the Lie algebra $L[S]$ is equal to $$ H_1(L[S])=
L[S]/([L[S],L[S]]).$$ Some of this structure is well-known, and
recorded next.
\begin{thm}\label{thm: homology of Lie algebras}
The abelianization of $L[x_1^{B}]$, $H_1L[x_1^{B}]$, is given by
the free abelian group with basis $x_1$, and
$[\cdots[x_1,x_{j_1}]x_{j_2}]\cdots ]x_{j_{q-1}}]x_{j_q}]$ for all
sequences $2 \leq j_1, j_2, \cdots, j_q$.

The abelianization of $L[\Lambda_n^{C}]$, $H_1L[\Lambda_n^{C}]$,
is given by the free abelian group with basis $\Lambda_n$, and
$[\cdots[\Lambda_n,\gamma_{j_1}]\gamma_{j_2}]\cdots ]
\gamma_{j_{q-1}}]\gamma_{j_q}]$ for $2 \leq j_1, j_2, \cdots,
j_q$.
\end{thm}

The graded free abelian groups  $H_1L[x_1^{B}]$, and
$H_1L[\Lambda_n^{C}]$ are filtered as described next for which the
filtration $H_1L[\Lambda_n^{C}]$ satisfies a reversed ordering
from that of $H_1L[x_1^{B}]$. This reversal is forced by the
infinitesimal braid relations, and the requirement that the map
$E_0^*(\Theta_n)$ preserve filtrations on the level of the first
homology group as suggested by Theorems \ref{thm:filtration of the
homology of a Lie algebras}, and \ref{thm: filtration preserving
on homology} below.

\begin{defn}\label{defn: filtrations for L...B}
Grade the first homology group $H_1(L[x_1^{B}])$ by setting
$$Gr_q(H_1(L[x_1^{B}])$$ for $ q \geq 1$ equal to the linear span
of the classes
$$[\cdots[x_1,x_{j_1}]x_{j_2}]\cdots ]x_{j_{q-1}}]x_{j_q}]$$
for all sequences $2 \leq j_1,j_2, \cdots j_{q-1},j_q $ with the
convention that $Gr_1(H_1(L[x_1^{B}])$ is the linear span of of
the class $x_1$. Thus
$$H_1(L[x_1^{B}]) = \oplus _{q \geq 1}Gr_q(H_1(L[x_1^{B}]).$$ Filter
$Gr_q(H_1(L[x_1^{B}])$ by $F_p(Gr_q(H_1(L[x_1^{B}]))$ as follows.
\begin{enumerate}
  \item For every $q \geq 1$,
  $$F_0Gr_q(H_1(L[x_1^{B}])$$ is the linear span of the classes
$$[\cdots[x_1,x_{j_1}]x_{j_2}]\cdots ]x_{j_{q-1}}]x_{j_q}]$$
for all sequences $$2\leq j_1\leq j_2\leq \cdots \leq j_q.$$
  \item Filtration $p > 0$, $F_p(Gr_q(H_1(L[x_1^{B}]))$, is
  the linear span of the classes
$$[\cdots[x_1,x_{j_1}]x_{j_2}]\cdots ]x_{j_{q}}]x_{j_{q}}]$$
with $D([\cdots[x_1,x_{j_1}]x_{j_2}]\cdots ]x_{j_{q}}]x_{j_{q}}])
\leq p$ where $$D([\cdots[x_1,x_{j_1}]x_{j_2}]\cdots
]x_{j_{q}}]x_{j_{q}}]) = \Sigma_{1 \leq i \leq q-1}d(x_{j_i})$$
with $d(x_{j_i})$ equal to the number of $x_{j_{i+k}}$ such that
$j_i > j_{i+k}$ for $ k> 0$.
\end{enumerate} Thus, there are inclusions
$$F_p(Gr_q(H_1(L[x_1^{B}])) \subset F_{p+1}(Gr_q(H_1(L[x_1^{B}]))$$
for all $p$, and $$Gr_q(H_1(L[x_1^{B}])) = \cup_{0 \leq p }
F_p(Gr_q(H_1(L[x_1^{B}])).$$
\end{defn}

There is a similar filtration $H_1L[\Lambda_n^{C}]$, but with
reversed ordering.
\begin{defn}\label{defn: filtration of Lie algebras of L...C}
Grade the first homology group $$H_1(L[\Lambda_n^{C}])$$ by
setting $$Gr_q(H_1(L[\Lambda_n^{C}])$$ for $ q \geq 1$ equal to
the linear span of the classes
$$[\cdots[\Lambda_n,\gamma_{j_1}]\gamma_{j_2}]\cdots ]
\gamma_{j_{q-1}}]\gamma_{j_q}]$$ for all sequences $2 \leq
j_1,j_2, \cdots j_{q-1},j_q $ with the convention that
$Gr_1(H_1(L[\Lambda_n^{C}])$ is the linear span of the class
$x_1$. Thus $$H_1(L[\Lambda_n^{C}]) = \oplus _{q \geq 1}
Gr_q(H_1(L[\Lambda_n^{C}]).$$ Filter $Gr_q(H_1(L[\Lambda_n^{C}])$
by $F_p(Gr_q(H_1(L[\Lambda_n^{C}]))$ as follows.
\begin{enumerate}
  \item For every $q \geq 1$,
  $$F_0Gr_q(H_1(L[\Lambda_n^{C}])$$ is the linear span of the classes
$$[\cdots[\Lambda_n,\gamma_{j_1}]\gamma_{j_2}]\cdots ]
\gamma_{j_{q-1}}]\gamma_{j_q}]$$ for all sequences $$ j_1\geq
j_2\geq \cdots \geq j_q \geq 2.$$
  \item Filtration $p > 0$, $F_p(Gr_q(H_1(L[\Lambda_n^{C}]))$ is
  the linear span of the classes
$$[\cdots[\Lambda_n,\gamma_{j_1}]\gamma_{j_2}]\cdots ]
\gamma_{j_{q-1}}]\gamma_{j_q}]$$ such that
$$\Delta([\cdots[\Lambda_n,\gamma_{j_1}]\gamma_{j_2}]\cdots ]
\gamma_{j_{q-1}}]\gamma_{j_q}])= \Sigma_{1 \leq i \leq q-1}
\delta(\gamma_{j_i})$$ with $\delta(\gamma_{j_i})$ equal to the
number of $\gamma_{j_{i+k}}$ such that $j_i < j_{i+k}$ for $ k>
0$.
\end{enumerate} Thus, there are inclusions
$$F_p(Gr_q(H_1(L[\Lambda_n^{C}])) \subset F_{p+1}(Gr_q(H_1(L[\Lambda_n^{C}]))$$
for all $p$, and
$$Gr_q(H_1(L[\Lambda_n^{C}])) = \cup_{0 \leq p }
F_p(Gr_q(H_1(L[\Lambda_n^{C}])).$$
\end{defn}

\begin{thm}\label{thm:filtration of the homology of a Lie algebras}
The map of Lie algebras $E_0^*(\Theta_n): L[x_1^{B}]\to\
L[\Lambda_n^{C}]$ induces a map $$H_1(\Theta_n):H_1 L[x_1^{B}]\to\
H_1L[\Lambda_n^{C}]$$ which
\begin{enumerate}
  \item preserves both gradation, and filtration as given in Definitions
\ref{defn: filtrations for L...B}, and \ref{defn: filtration of
Lie algebras of L...C}, and
  \item is an isomorphism.
\end{enumerate}
\end{thm}

\begin{cor}\label{cor:Theta is an iso on associated graded}
The map of Lie algebras $E_0^*(\Theta_n): L[x_1^{B}]\to\
L[\Lambda_n^{C}]$ is an isomorphism.
\end{cor}

The proof of Theorem \ref{thm:filtration of the homology of a Lie
algebras} follows by induction with the first case given next in
which $x_i = x_i(n)$, and $\gamma_i = \gamma_i(n)$.
\begin{thm}\label{thm: admissible sequences}
If $n \geq 1$, then $E_0^*(\Theta_n)(x_1) = \Lambda_n$. If $2 \leq
j_1 \leq j_2 \leq \cdots \leq j_q$, then
$$E_0^*(\Theta_n)([\cdots[x_1,x_{j_1}]x_{j_2}]\cdots ]
x_{j_{q-1}}]x_{j_q}]) = [\cdots[\Lambda_n
,\gamma_{j_q}]\gamma_{j_{q-1}}] \cdots ]
\gamma_{j_{2}}]\gamma_{j_1}].$$ Thus $E_0^*(\Theta_n)$ preserves
filtration $0$ in Theorem \ref{thm:filtration of the homology of a
Lie algebras}, and induces an isomorphism
$$H_1(\Theta_n):F_0(Gr_q(H_1(L[x_1^{B}]))  \to\
F_0(Gr_q(H_1(L[\Lambda_1^{C}]))$$ for all $ q \geq 0$.
\end{thm}

\begin{proof}
By Theorem \ref{thm: values of theta}, $E_0^*(\Theta_n)(x_1) =
\Lambda_n$. Thus the theorem is correct in the case of the empty
sequence. The next step is to check that $E_0^*(\Theta_n)$
preserves filtration $0$ by inducting on $q$ as stated in Theorem
\ref{thm:filtration of the homology of a Lie algebras}. Assume by
induction that if $2 \leq j_1 \leq j_2 \leq \cdots \leq j_q$, then
$$E_0^*(\Theta_n)([\cdots[x_1,x_{j_1}]x_{j_2}]\cdots ]
x_{j_{q-1}}]x_{j_q}]) = [\cdots[\Lambda_n
,\gamma_{j_q}]\gamma_{j_{q-1}}] \cdots ]
\gamma_{j_{2}}]\gamma_{j_1}].$$

Next, consider $j_q \leq j_{q+1}$, together with the value of
$$E_0^*(\Theta_n)([\cdots[x_1,x_{j_1}]x_{j_2}]\cdots ]
x_{j_{q-1}}]x_{j_q}] x_{j_{q+1}}])$$ given by
$$[\cdots[\Lambda_n ,\gamma_{j_q}]\gamma_{j_{q-1}}]
\cdots ] \gamma_{j_{2}}]\gamma_{j_1}]E_0^*(\Theta_n)(x_{q+1})].$$
If $2 \leq j_p \leq j_q$, then $[\gamma_{j_p}(n),
\Theta(x_{j_q}(n))]= 0$ by Theorem \ref{thm: values of theta}.
Hence by the Jacobi identity, $[[A,B],C] = [[A,C],B] + [A,[B,C]]$,
the value of $$[\cdots[\Lambda_n ,\gamma_{j_q}]\gamma_{j_{q-1}}]
\cdots ] \gamma_{j_{2}}]\gamma_{j_1}]E_0^*(\Theta_n)(x_{q+1})]$$
is equal to
$$[\cdots[\Lambda_n ,\gamma_{j_q}]\gamma_{j_{q-1}}]
\cdots ] \gamma_{j_{2}}]E_0^*(\Theta_n)(x_{q+1})]\gamma_{j_1}].$$

Notice that by the inductive hypothesis, the element
$$[\cdots[\Lambda_n ,\gamma_{j_q}]\gamma_{j_{q-1}}]
\cdots ] \gamma_{j_{2}}]E_0^*(\Theta_n)(x_{q+1})]$$ is equal to
$$E_0^*(\Theta_n)([\cdots[x_1,x_{j_2}]x_{j_3}]\cdots ]
x_{j_{q-1}}] x_{j_{q+1}}]).$$ That is, the element $x_{j_1}$ does
not appear.

Furthermore, notice that the element
$E_0^*(\Theta_n)([\cdots[x_1,x_{j_2}]x_{j_3}]\cdots ] x_{j_{q-1}}]
x_{j_{q+1}}])$ is equal to
$$[\cdots[\Lambda_n ,\gamma_{j_{q+1}}]\gamma_{j_{q-1}}]
\cdots ] \gamma_{j_{3}}]\gamma_{j_2}]$$ by the inductive
hypothesis. Hence
$$E_0^*(\Theta_n)([\cdots[x_1,x_{j_1}]x_{j_2}]\cdots ]
x_{j_{q-1}}]x_{j_q}] x_{j_{q+1}}])= [\cdots[\Lambda_n
,\gamma_{j_{q+1}}]\gamma_{j_{q}}] \cdots ]
\gamma_{j_{2}}]\gamma_{j_1}].$$ The theorem follows.
\end{proof}

The proof of the next theorem is analogous, but uses a variation
of Theorem \ref{thm: values of degeneracies and theta} to measure
changes of order, and their effect on the filtration defined
above. To keep track of certain signs, the following conventions
are used.

\begin{thm}\label{thm: filtration preserving on homology}
If $2 \leq j_1, j_2, \cdots, j_q$, then the class of
$$E_0^*(\Theta_n)([\cdots[x_1,x_{j_1}]x_{j_2}]\cdots ]
x_{j_{q-1}}]x_{j_q}])$$ in $H_1L[\Lambda_n^{C}]$ is equal to the
class of $$ (\pm 1) [\cdots[\Lambda_n
,\gamma_{j_q}]\gamma_{j_{q-1}}]\cdots ]
\gamma_{j_{2}}]\gamma_{j_1}] + \Omega$$ where $\Omega$ projects to
an element of lower filtration degree in $H_1L[\Lambda_n^{C}]$. In
addition, $E^*_0(\Theta_n)$ both preserves the filtrations as
given in Definitions \ref{defn: filtrations for L...B}, and
\ref{defn: filtration of Lie algebras of L...C} , as well induces
an isomorphism on the level of first homology groups.
\end{thm}

\begin{proof}
The proof of Theorem \ref{thm: filtration preserving on homology}
is by induction on filtration degree $p$ for fixed gradation
$Gr_q(H_1(L[x_1^{B}])$ starting with filtration degree $p = 0$.
Notice that if $n \geq 1$, then $E_0^*(\Theta_n)(x_1) =
\Lambda_n$, and if $2 \leq j_1 \leq j_2 \leq \cdots \leq j_q$,
then
$$E_0^*(\Theta_n)([\cdots[x_1,x_{j_1}]x_{j_2}]\cdots ]
x_{j_{q-1}}]x_{j_q}]) = [\cdots[\Lambda_n
,\gamma_{j_q}]\gamma_{j_{q-1}}] \cdots ]
\gamma_{j_{2}}]\gamma_{j_1}]$$ by \ref{thm: admissible sequences}.
Thus the theorem is correct for filtration degree $0$, and for all
$q$ as the map $E_0^*(\Theta_n)$ induces an isomorphism on the
level of $$H_1(\Theta_n):F_0Gr_q(H_1(L[x_1^{B}])) \to\
F_0Gr_q(H_1(L[\Lambda_n^{C}]))$$ for all $1 \leq q$.

Thus assume inductively that
$$E_0^*(\Theta_n)[\cdots[x_1,x_{j_1}]x_{j_2}]\cdots ]
x_{j_{q-1}}]x_{j_q}]$$ in $H_1L[\Lambda_n^{C}]$ is equal to $(\pm
1 ) [\cdots[\Lambda_n ,\gamma_{j_q}]\gamma_{j_{q-1}}]\cdots ]
\gamma_{j_{2}}]\gamma_{j_1}] + \Omega$ where $\Omega$ is of lower
filtration as stated in Theorem \ref{thm: filtration preserving on
homology}.

Consider the following formulae.
\begin{enumerate}
  \item $E_0^*(\Theta_n)[\cdots[x_1,x_{j_1}]x_{j_2}]\cdots ]
x_{j_{q-1}}]x_{j_q}]x_{j_{q+1}}]  = [[A,B],C]$ for which
$$A = E_0^*(\Theta_n)([\cdots[x_1,x_{j_1}]x_{j_2}]\cdots ]
x_{j_{q-1}}])$$ with $$B=E_0^*(\Theta_n)(x_{j_{q}}),$$ and
$$C= E_0^*(\Theta_n)(x_{j_{q+1}}).$$
  \item The inductive hypothesis gives that
$$[A,B] = (\pm 1)[\cdots[\Lambda_n
,\gamma_{j_q}]\gamma_{j_{q-1}}]\cdots ]
\gamma_{j_{2}}]\gamma_{j_1}] + \Omega $$ for which $\Omega$
projects to lower filtration degree in $Gr_q(H_1(L[x_1^{B}))$.
Hence $$[[A,B],C] = [((\pm 1)[\cdots[\Lambda_n
,\gamma_{j_q}]\gamma_{j_{q-1}}]\cdots ]
\gamma_{j_{2}}]\gamma_{j_1}] + \Omega)],C]$$ in
$Gr_q(H_1(L[x_1^{C}])$.

\item Since $\Omega$ has lower filtration degree than $[A,B]$, the
filtration degree of $[\Omega,C]$ has filtration degree strictly
less than that of $[[A,B],C]$ in $Gr_q(H_1(L[\Lambda_n^{C}])$ by
inspection of the definition.

\item Next consider $$[\cdots[\Lambda_n ,\gamma_{j_q}]
\gamma_{j_{q-1}}]\cdots ] \gamma_{j_{2}}]\gamma_{j_1}],C] =
[[E,F]C]= [[E,C]F] + [E[F,C]]$$ for which
\begin{itemize}
  \item $E = [\cdots[\Lambda_n ,\gamma_{j_q}]\gamma_{j_{q-1}}]\cdots ]
\gamma_{j_{2}}]$, and
  \item $F = \gamma_{j_{q-1}}$.
\end{itemize}

\item Notice that

$$[E,C] = [\cdots[\Lambda_n ,\gamma_{j_q}]\gamma_{j_{q-1}}]\cdots ]
\gamma_{j_{2}}]E_0^*(\Theta_n)(x_{j_{q+1}})].$$

\item Thus $[E,C]$ is the sum
$$[[X,E_0^*(\Theta_n)(x_{j_{q+1}})]\gamma_{j_{2}}] +
[X,[\gamma_{j_{2}},E_0^*(\Theta_n)(x_{j_{q+1}})]]$$ where $$X =
\cdots[\Lambda_n ,\gamma_{j_q}]\gamma_{j_{q-1}}]\cdots ]
\gamma_{j_{3}}].$$

\item If $$j_2 \leq j_{q+1},$$ then
$[\gamma_{j_{2}},E_0^*(\Theta_n)(x_{j_{q+1}})] = 0$ by part $6$ of
Theorem \ref{thm: values of degeneracies and theta}. Hence
$$[E,C]= [[X,E_0^*(\Theta_n)(x_{j_{q+1}})]\gamma_{j_{2}}],$$ the
inductive hypothesis applies, and the Theorem follows.

\item  If  $$j_2 > j_{q+1},$$ then
$$[\gamma_{j_{2}},E_0^*(\Theta_n)(x_{j_{q+1}})] = V + Y$$ where $V
= [\gamma_j(2),\Lambda_n]+ [\Lambda_n,\gamma_{j_{q+1}}]$, and $Y =
2[\gamma_{j_{2}},\gamma_{j_{q+1}}]$ by part $5$ of Theorem
\ref{thm: values of degeneracies and theta}. Thus
$[[X,E_0^*(\Theta_n)(x_{j_{q+1}})],\gamma_{j_{2}}]
+[X,[\gamma_{j_{2}},E_0^*(\Theta_n)(x_{j_{q+1}})]]$ is equal to
the coset of $$[[X,E_0^*(\Theta_n)(x_{j_{q+1}})]\gamma_{j_{2}}] +
[X,2[\gamma_{j_{2}},\gamma_{j_{q+1}})]]$$ in $H_1(L[x_1^{C}])$,
and thus  $$[[X,E_0^*(\Theta_n)(x_{j_{q+1}}) -
2\gamma_{j_{q+1}}]\gamma_{j_{2}}]$$ in $H_1(L[x_1^{C}])$ modulo
terms of lower filtration.

\item Furthermore, $$[[X,E_0^*(\Theta_n)(x_{j_{q+1}}) -
2\gamma_{j_{q+1}}]\gamma_{j_{2}}] =
-[[X,\gamma_{j_{q+1}}]\gamma_{j_{2}}]$$ in case $X = \Lambda_n$,
and

$$[[\gamma_s,E_0^*(\Theta_n)(x_{j_{q+1}}) -
2\gamma_{j_{q+1}}]\gamma_{j_{2}}] = 0.$$

\item Thus if $j_2 > j_{q+1},$ then
$[[X,E_0^*(\Theta_n)(x_{j_{q+1}}) -
2\gamma_{j_{q+1}}]\gamma_{j_{2}}]$ is one of the following.
        \begin{enumerate}
        \item -$[[X,\gamma_{j_{q+1}}]\gamma_{j_{2}}]$, or
        \item $[[X,\gamma_{j_{q+1}}]\gamma_{j_{2}}]$
        \end{enumerate}
\item In addition, $$[\Lambda_n,E_0^*(\Theta_n)(x_{j_{q+1}})]
=[\Lambda_n,\gamma_{j_{q+1}}]$$ as follows from \ref{thm: values
of degeneracies and theta}, and \ref{thm: values of theta}.
\end{enumerate}

Formula $5$ in Theorem \ref{thm: values of degeneracies and theta}
is used at this point, and is stated next for the convenience of
the reader. For any fixed $2 \leq i< j \leq n$,
$[\gamma_j(n),E_0^*(\Theta_n)(x_i)] = S$ where $S =
[\gamma_j(n),\Lambda_n]+ [\Lambda_n,\gamma_i(n)] +
 2[\gamma_j(n),\gamma_i(n)].$ This formula will be used in the
expansion of $E_0^*(\Theta_n)[\cdots[x_1,x_{j_1}]x_{j_2}]\cdots ]
x_{j_{q-1}}]x_{j_q}]x_{j_{q+1}}]$. Thus,
$$[B,C]= [E_0^*(\Theta_n)(x_{j_{q}}),E_0^*(\Theta_n)(x_{j_{q+1}})]
= [\gamma_{j_q},\Lambda_n]+ [\Lambda_n,\gamma_{j_{q+1}}] +
 2[\gamma_{j_{q}},\gamma_{j_{q+1}}].$$

Formula $5$ in Theorem \ref{thm: values of degeneracies and theta}
is used at this point, and is stated next for the convenience of
the reader. For fixed $i$, and $j$ with $2 \leq i< j \leq n$,
$$[\gamma_j(n),E_0^*(\Theta_n)(x_i)] = U + V$$
where $U = [\gamma_j(n),\Lambda_n]+ [\Lambda_n,\gamma_i(n)]$, and
$V = 2[\gamma_j(n),\gamma_i(n)]$.

The previous formula will be used in the expansion of the element
specified by an expansion of $E_0^*(\Theta_n)([X_q,x_{j_{q+1}}])$
for $X_q = [\cdots[x_1,x_{j_1}]x_{j_2}]\cdots
]x_{j_{q-1}}]x_{j_q}]$. Thus, $[B,C]=
[E_0^*(\Theta_n)(x_{j_{q}}),E_0^*(\Theta_n)(x_{j_{q+1}})] = U"
+V"$ where $U" = [\gamma_{j_q},\Lambda_n]+
[\Lambda_n,\gamma_{j_{q+1}}]$, and $V" =
2[\gamma_{j_{q}},\gamma_{j_{q+1}}]$.

Next, notice that $$[A,[\gamma_q,\Lambda_n]],$$ and
$$[A,[\Lambda_n, \gamma_{q+1}]]$$ both project to $0$ in
$Gr_q(H_1(L[\Lambda_n^{C}])$. Hence the class of $[A,[B,C]]$ equal
to the class of $[A,2[\gamma_{j_{q}},\gamma_{j_{q+1}}]]$.

\end{proof}

\begin{thm}\label{thm: isomorphism}
The map
$$ E_0^*(\Theta_n): L[x_1^{B}]\to\
L[\Lambda_n^{C}]$$ induces an isomorphism of Lie algebras, and
thus a monomorphism. Hence $\Theta_n$ is a monomorphism.
\end{thm}

\begin{proof}
Notice that the map $E_0^*(\Theta_n)$ sends a generator
$[\cdots[x_1,x_{j_1}]x_{j_2}]\cdots ] x_{j_{q-1}}]x_{j_q}]$ to
$$(\pm 1) [\cdots[\Lambda_n
,\gamma_{j_q}]\gamma_{j_{q-1}}]\cdots ]
\gamma_{j_{2}}]\gamma_{j_1}] + \Omega$$ where $\Omega$ projects to
an element of lower filtration degree in $H_1L[\Lambda_n^{C}]$ by
Theorem \ref{thm: filtration preserving on homology}. Hence the
morphism of Lie algebras $E_0^*(\Theta_n)$ induces an isomorphism
on the module of indecomposables, and thus an isomorphism of Lie
algebras. The Theorem follows.
\end{proof}

\section{On Vassiliev invariants, the mod-$p$ descending central
series, and the Bousfield-Kan spectral sequence} \label{Vassiliev
invariants, and the Bousfield-Kan spectral sequence}

Let $\Gamma^n(G)$, respectively $\Gamma^{n,p}(G)$, denote the
$n$-th stage of the descending central series for a discrete group
$G$, respectively, the mod-$p$ descending central series for $G$.
Thus
\begin{enumerate}
    \item $\Gamma^n(G)$ is the subgroup of $G$ generated by commutators
$[\cdots[g_1,g_2],g_3], \cdots],g_t]$ for $ t \geq n$ with
decreasing filtration
$$G = \Gamma^1(G) \supseteq \Gamma^2(G) \supseteq \cdots, $$
 and
    \item $\Gamma^{n,p}(G)$ is the subgroup of $G$ generated by commutators
$[\cdots[g_1,g_2],g_3], \cdots],g_s]^{p^j}$ for $ s\cdot p^j \geq
n$ with decreasing filtration
$$G = \Gamma^{1,p}(G) \supseteq \Gamma^{2,p}(G) \supseteq \cdots .$$
\end{enumerate}

Let $E_0^n(G) = \Gamma^n(G)/ \Gamma^{n+1}(G)$, and $E_{0}^{n,p}(G)
= \Gamma^{n,p}(G)/ \Gamma^{n+1,p}(G)$. The commutator map of sets
$$[-,-]: G \times G \to\ G$$ induces a natural pairing endowing
$$E_0^*(G) = \oplus_{n \geq 1} E_0^n(G) $$ with
the structure of a Lie algebra while the analogous pairing for
$$E_{0}^{*,p}(G) = \oplus_{n \geq 1}E_{0}^{n,p}(G) $$
with the $p$-th power map $\xi: G \to G$ which induces a function
$$\xi:E_{0}^{n,p}(G) \to E_{0}^{pn,p}(G)$$ gives the structure of a restricted
Lie algebra \cite{jac}.

Results of Kohno-Falk-Randell used above give the structure of
$E_0^n(P_{n+1})$ as well as a relationship between the Lie algebra
obtained from the descending central series of a group with
certain natural choices of group extensions. The variation given
in this section replaces the descending central series by the
mod-$p$ descending central series thus giving the analogous
structure for $E_{0}^{*,p}(P_{n+1})$.

Recall that the unstable Adams spectral sequence is that obtained
by filtering a simplicial group by the mod-$p$ descending central
series \cite{Cu}. This variation is recorded here as it
corresponds to the relationship between Vassiliev invariants, the
Bousfield-Kan spectral sequence, or unstable Adams spectral
sequence.

First, consider the integral version gotten by filtering by the
descending central series which is labelled the Bousfield-Kan
spectral sequence below. On the level of $E_0^*$, the morphism of
simplicial groups $\Theta: F[S^1] \to\ \mathrm{AP}_*$ induces a
map
$$E_0^*(\Theta_n): E_0^*(F_{n}) \to\ E_0^*(P_{n+1}),$$ the subject
of sections \ref{sec: lie algebras}, and \ref{sec:On Theta n} on
embeddings of Lie algebras here. In addition, the associated
graded Lie algebra for the mod-$p$ descending central series
$E_0^{*,p}(F_{n})$ gives the $E_0$-term of a spectral sequence
abutting to the homotopy groups of $F[S^1]$ modulo torsion prime
to $p$.

On the other hand, the Lie algebra $E_0^*(P_{n+1})$ was studied in
\cite{K,K1,K2} where the universal enveloping algebra was shown to
give Vassiliev invariants of pure braids. The Vassiliev invariants
distinguish pure braids as proven in \cite{K}. Thus they
distinguish elements in the $E^0$-term of the Bousfield-Kan
spectral sequence for $F[S^1]$. The analogue for the mod-$p$
descending central series is developed next.

\noindent{\bf Question:} Is there a further relationship between
the Vassiliev invariants of pure braids, and the homotopy groups
of the $2$-sphere ? Is there an informative interplay between
these invariants, and homotopy theory ?

The point of the next remarks is to record the structure of the
Lie algebra obtained from the mod-$p$ descending central series
for the pure braid groups. Algebraic preliminaries are given next
arising from work of \cite{K,FR,Xico}, and the modification below
for the mod-$p$ descending central series.

\begin{thm}\label{thm:split short exact sequence of groups and Lie algebras}
Let
\[
\begin{CD}
1 \to\ A @>{j}>> B @>{p}>> C \to\ 1
\end{CD}
\] be a short exact sequence of groups such that
\begin{itemize}
\item there is a section $\sigma$ for $p: B  \to\ C$ giving $p
\circ \sigma = 1_C$, and \item the natural action of C on $H_1(A)
$ is trivial: Given $b$ in $B$ and $a$ in $A$, then $$bab^{-1} =
ax$$ for some $x$ in the commutator subgroup $[A,A] =
\Gamma^2(A)$.
\end{itemize}

Then there is

\begin{enumerate}
 \item a split short exact sequences of Lie algebras
 $$ 0 \to\ E_0^*(A) \to\  E_0^*(B) \to\  E_0^*(C) \to\ 0,$$ and
 \item a split short exact sequences of restricted Lie algebras
$$ 0 \to\ E_{0}^{*,p}(A) \to\  E_{0}^{*,p}(B) \to\ E_{0}^{*,p}(C) \to\ 0.$$
\end{enumerate}
\end{thm}

The structure of the Lie algebra for the mod-$p$ descending
central series of the pure braid group, as well as certain other
groups follows from the proof of the Proposition \ref{prop: Lie
algebras}. A proof is analogous to the ones in \cite{K,FR,Xico};
modifications in the case of the mod-$p$ descending central series
are direct, and are listed below for convenience.

Observe that $p((b\cdot \sigma (p(b^{-1}))) = 1$, and so there
exists an unique element $a$ in A with $j(a) = b\cdot \sigma
p(b^{-1})$. Thus, there is a well-defined function ( not
necessarily a homomorphism )
$$ \tau: B \to\ A $$ defined by
the formula $$ \tau(b) = j^{-1}(b \cdot \sigma (p(b^{-1}))).$$ A
useful variation of \ref{thm:split short exact sequence of groups
and Lie algebras} for the mod-$p$ descending central series is
recorded next.

\begin{prop}\label{prop: Lie algebras}
Let
\[
\begin{CD}
1 \to\ A @>{j}>> B @>{p}>> C \to\ 1
\end{CD}
\] be a short exact sequence of groups such that
\begin{itemize}
\item there is a section $\sigma$ for $p: B  \to\ C$ giving $p
\circ \sigma = 1_C$, and \item the natural action of C on $H_1(A)
$ is trivial: Given $b$ in $B$ and $a$ in $A$, then $$bab^{-1} =
ax$$ for some $x$ in the commutator subgroup $[A,A]= \Gamma^2(A)$.
( Note that $x$ is always in $[B,B]$, but $x$ is not necessarily
an element in $[A,A]$.)
\end{itemize}

Then the following hold.
\begin{enumerate}
\item  For fixed $b$ in $B$, $a$ in $A$, there is an element $y$
in $[A,A]$ such that $b^{-1}ab = ay$.

\item  The group $[B,A]$ is a subgroup of $[A,A]$.

\item  The group $[\Gamma^{m}(B),\Gamma^{n}(A)]$ is a subgroup of
$\Gamma^{n+m}(A)$, and the group
$[\Gamma^{m,p}(B),\Gamma^{n,p}(A)]$ is a subgroup of
$\Gamma^{n+m,p}(A)$.

\item  If $ n \geq 1 $, $\tau$ is a filtration preserving function
( not necessarily a group homomorphism ). That is,
$\tau(\Gamma^{n}(B))$ is contained in $\Gamma^{n}(A)$, and
$\tau(\Gamma^{n,p}(B))$ is contained in $\Gamma^{n,p}(A)$.

\item If $p(b) = 1$, then $j\tau(b)  = b $. Furthermore, there are
split short exact sequences of groups

$$1 \to\ \Gamma^n(A) \to\ \Gamma^n(B) \to\ \Gamma^n(C) \to\ 1,$$ and
$$1 \to\ \Gamma^{n,p}(A) \to\ \Gamma^{n,p}(B) \to\ \Gamma^{n,p}(C) \to\ 1.$$

\item If $ n \geq 1 $, there are well-defined induced functions
$\overline\tau: E_0^n(B) \to\ E_0^n(A),$ and $\overline\tau:
E_0^{n,p}(B) \to\ E_0^{n,p}(A)$ defined on an equivalence class of
$b$, $[b]$, by the formula $\overline\tau([b])= [\tau(b)]$.

\item If $ n \geq 1 $, and $[b]$ is in the kernel of the induced
homomorphism $E_0^n(p):E_0^n(B) \to\ E_0^n(C),$ respectively $[b]$
is in the kernel of the induced homomorphism
$E_0^{n,p}(p):E_0^{n,p}(B) \to\ E_0^{n,p}(C),$ then
$E_0^n(j)(\overline\tau [b])$ = $[b]$, respectively
$E_0^{n,p}(j)(\overline\tau [b])$ = $[b]$ .

Furthermore, there are split short exact sequences of abelian
groups
$$ 0 \to\ E_0^n(A) \to\  E_0^n(B) \to\  E_0^n(C) \to\ 0,$$ and
$$ 0 \to\ E_{0}^{n,p}(A) \to\  E_{0}^{n,p}(B) \to\  E_{0}^{n,p}(C) \to\ 0.$$
\item The morphisms $ 0 \to\ E_{0}^{*,p}(A) \to\  E_{0}^{*,p}(B)
\to\ E_{0}^{*,p}(C) \to\ 0$ are of restricted Lie algebras.
\end{enumerate}
\end{prop}

\begin{remark}
The hypotheses above that $bab^{-1}$ = $ax$ for some $x$ in the
commutator subgroup $[A,A]$ is not necessarily satisfied without
the hypotheses of trivial local coefficients. Notice that if $A$
is a normal subgroup of $B$, then $a^{-1}b^{-1}ab$ is always in
$A$, but may not necessarily be in $[A,A]$. An example is given by
the group extension $ 1 \to\ Z/3Z \to\  \Sigma_3  \to\ Z/2Z  \to\
1$ for which the commutator subgroup $A$ = $Z/3Z$ is trivial, but
the group $[B,A]$ is non-trivial for $B$ = $\Sigma_3$. The
hypothesis of trivial local coefficients is important in this
step.
\end{remark}

Theorem \ref{thm:split short exact sequence of groups and Lie
algebras} is a restatement of Proposition \ref{prop: Lie
algebras}. The proof of \ref{prop: Lie algebras} is based on the
Hall-Witt identities together with another lemma of P.~Hall both
recorded in the next statement proven in \cite{MKS}, page $290$,
and \cite{DDMS}, page $2$. An additional useful statement is also
given by \ref{thm:Hall lemma}. Recall that $[a,b]= a^{-1} \cdot
b^{-1} \cdot a \cdot b$ denotes the commutator of elements $a$,
and $b$ in a group $G$, and $c^a = c^{-1}\cdot a \cdot c$.

\begin{thm}\label{thm:Hall lemma}
For any elements of a group $G$,

\begin{enumerate}
    \item $[a,b]\cdot [b,a] = 1$,

    \item $[a,b\cdot c] = [a,c]\cdot [a,b]\cdot [[a,b],c]$,
    \item $[a \cdot b, c] = [a,c]\cdot [[a,c],b]\cdot [b,c]$,

    \item $[[a, b], c^a] \cdot [[c, a], b^c] \cdot [[b, c], a^b] =
    1$,

    \item $[[a, b],c] \cdot [[b,c],a] \cdot [[c,a],b] =
    [b,a]\cdot [c,a]\cdot [c,b]^a \cdot [a,b]\cdot [a,c]^b \cdot [b,c]^a \cdot
    [a,c]\cdot [c,a]^b$, and

    \item $[a,b^n] = [a,b]\cdot [a,b]^b \cdots [a,b]^{b^{n-1}}$.

If $a$ is an element of $\Gamma^{m,p}(G)$, and $b$ is an element
of $\Gamma^{n,p}(G)$, then  $$[a, b^p] = [a,b]^p\cdot z$$ for $z$
in $\Gamma^{(m+pn),p}(G)$.
\end{enumerate}
If $A$, $B$, and $C$ are normal subgroups of a group $G$, then
$[[A,B],C]$ is contained in the subgroup generated by $[[B,C],A]$
and $[[C,A],B]$.
\end{thm}

The proof of Proposition \ref{prop: Lie algebras} is given next.
\begin{proof} The first statement of the
proposition is one of the stated assumptions concerning the
extension $ 1 \to\ A  \to\   B  \to\ C \to\ 1 $. That is $b^{-1}ab
= ax$ for $b$ in B, a in $A$, and some $x$ in $[A,A]$ by the
assumption that the local coefficient system on $H_1(A)$ is
trivial.

By part (1), $b^{-1}ab = ax$ for $b$ in $B$, $a$ in $A$, and some
$x$ in $[A,A]$. Thus $a^{-1}b^{-1}ab = x $, so $[B,A]$ is a
subgroup of $[A,A]$. Statement (2) follows.

The proof that $[\Gamma^{m}(B),\Gamma^{n}(A)]$ is a subgroup of
$\Gamma^{n+m}(A)$ is given in \cite{FR,Xico}. Modifications in the
case of $[\Gamma^{m,p}(B),\Gamma^{n,p}(A)]$ are listed next. Since
$\Gamma^{1}(G)= \Gamma^{1,p}(G)$ for any group $G$,
$[\Gamma^{1,p}(B),\Gamma^{1,p}(A)] =[\Gamma^{1}(B),\Gamma^{1}(A)]$
is a subgroup of $\Gamma^{2}(A)$, and thus $\Gamma^{2,p}(A)$.
Consider the case of $[\Gamma^{m,p}(B),\Gamma^{n,p}(A)]$.
Inductively assume that for $m< M$, and $n<N$ that
$[\Gamma^{m,p}(B),\Gamma^{n,p}(A)]$ is a subgroup of
$\Gamma^{n+m,p}(A)$. That $[\Gamma^{m,p}(B),\Gamma^{n+q,p}(A)]$ is
also a subgroup of $\Gamma^{n+m+q,p}(A)$ follows by induction on
$q$ via Theorem \ref{thm:Hall lemma}, or Appendix A of
\cite{DDMS}: If $a$ is an element of $\Gamma^{m,p}(G)$, and $b$ is
an element of $\Gamma^{n,p}(G)$, then  $[a, b^p] = [a,b]^p\cdot z$
for $z$ in $\Gamma^{(m+pn),p}(G)$ by Theorem \ref{thm:Hall lemma}.
Assume that $b$ is an element of $\Gamma^{m,p}(B)$, and that $a^p$
is an element of $\Gamma^{n+q-1,p}(A)$, then $[b,a^p]$ is an
element of $\Gamma^{m+n+q,p}(A)$.

The other cases are those in \cite{K,FR,Xico}: If $G =
\Gamma^{m,p}(B)$, $K = A$, $H= \Gamma^{n+q-1,p}(A)$, then
$[\Gamma^{m,p}(B),\Gamma^{n+q,p}(A)]$ is a subgroup of
$\Gamma^{n+m+q,p}(A)$. A similar argument implies that if
$[\Gamma^{m,p}(B),\Gamma^{n,p}(A)]$ is a subgroup of
$\Gamma^{n+m,p}(A)$ for $m <M$, and $n < N$, then
$[\Gamma^{m+s,p}(B),\Gamma^{n,p}(A)]$ is a subgroup of
$\Gamma^{n+m+s,p}(A)$.

Recall that there is a well-defined function $ \tau: B \to\ A $
given by the formula $ \tau(b) = j^{-1}(b \cdot \sigma
(p(b^{-1}))).$ To prove the fourth statement, notice that the
following formula holds which measures the failure of the map
$\tau$ from being multiplicative:
$$ \tau(xy) =  (\tau(x))(\tau(y))j^{-1}([(\tau(y))^{-1},\sigma
p(x)]).$$ In addition, $[(\tau(y))^{-1},\sigma p(x)])$ is an
element of $[A,B]$ which is a subgroup of $[A,A]$. Thus for all
elements $x,y$ in $B$, $$ \tau(xy) = (\tau(x))(\tau(y))V,$$ as
well as $$ \tau(x^p) = (\tau(x))^pW$$ where $V$, and $W$ are
elements of $[A,A]$ by statement (3).

Let $\lambda = \sigma \circ p$,  and observe that
$$ \tau(xyx^{-1}y^{-1}) = j^{-1}((xyx^{-1}y^{-1})(\lambda(y^{-1})(\lambda
x^{-1})(\lambda x)(\lambda y)).$$ Furthermore, if $v$ is an
element of $\Gamma^{m}(B)$, and $\tau(z)$ lies in $\Gamma^{n}(A)$,
then the commutator $[\tau(z),v]$ lies in $\Gamma^{m+n}(A)$ by
part (3). Similarly, if $v$ is an element of $\Gamma^{m,p}(B)$,
and $\tau(z)$ lies in $\Gamma^{n,p}(A)$, then the commutator
$[\tau(z),v]$ lies in $\Gamma^{m+n,p}(A)$ by part (3). Thus the
statement that $\tau(\Gamma^{n}(B))$ is contained in
$\Gamma^{n}(A)$, and $\tau(\Gamma^{n,p}(B))$ is contained in
$\Gamma^{n,p}(A)$ follows by induction on $n$ together with part
(3) of the proposition. Statement (4) follows.

Statement (5) concerns $ \tau: B \to\ A $ given by the formula $
\tau(b) = j^{-1}(b \cdot \sigma (p(b^{-1}))).$ By definition, if
$p(b) = 1$, then $p(b^{-1}) = 1$, and $j\tau(b)  = b$. Hence if if
$b$ is either in $\mathrm{ker}(p)\cap \Gamma^n(B)$, or
$\mathrm{ker}(p)\cap \Gamma^{n,p}(B)$, then $\tau(b) = b$. Thus by
the preceding remark as well as part (4), $ \tau$ restricts to
functions $ \tau|_{\mathrm{ker}(p)\cap \Gamma^n(B)
}:\mathrm{ker}(p)\cap \Gamma^n(B) \to\ \Gamma^n(A)$, and $
\tau|_{\mathrm{ker}(p)\cap \Gamma^{n,p}(B) }:\mathrm{ker}(p)\cap
\Gamma^{n,p}(B) \to\ \Gamma^{n,p}(A)$.

It follows that the homomorphisms $$ j: \Gamma^n(A) \to\
\mathrm{ker}(p)\cap \Gamma^n(B),$$ as well as $$ j:
\Gamma^{n,p}(A) \to\ \mathrm{ker}(p)\cap \Gamma^{n,p}(B)$$ are
group isomorphisms. Furthermore, there are exact sequences of
groups $ 1 \to\ \Gamma^n(A) \to\ \Gamma^n(B) \to\ \Gamma^n(C) \to\
1$, and $ 1 \to\ \Gamma^{n,p}(A) \to\ \Gamma^{n,p}(B) \to\
\Gamma^{n,p}(C) \to\ 1$ which are split by the existence of
$\sigma$. Part (5) follows.

If $ n \geq 1 $, there is a well-defined induced function
$\overline\tau: E_0^n(B) \to\ E_0^n(A)$ defined on an equivalence
class of b by the formula $\overline\tau([b])= \tau(b)$ by parts
(c), and (d) together with the formula $$ \tau(xy)
(\tau(x))(\tau(y))j^{-1}([(\tau(y))^{-1},\sigma p(x)]).$$

If the class $[b]$ of an element $b$ in $\Gamma^{n,p}(B)$ is in
the kernel of the induced homomorphism $E_0^{n,p}(p):E_0^{n,p}(B)
\to\ E_0^{n,p}(C)$, then $b$ is in $\Gamma^{n,p}(B)$, and $p(b)$
is in $\Gamma^{n+1,p}(C)$.

Since the natural map $\Gamma^{n+1,p}(p): \Gamma^{n+1}(B) \to\
\Gamma^{n+1}(C)$ is a surjection by part (3), there is an element
$x$ in $\Gamma^{n+1,p}(B)$ such that $p(x)$ = $p(b)$. Note that
$b$ is in $\Gamma^{n,p}(B)$, but $x$ is in $\Gamma^{n+1,p}(B)$.

Hence $x^{-1}$ is in $\Gamma^{n+1,p}(B)$, and $bx^{-1}$ is a
representative of the class $[b]$ in $E_0^n(B)$ with $j\tau(bx) =
bx$. Hence the natural map $E_0^{n,p}(j):E_0^{n,p}(A) \to\
E_0^{n,p}(B)$ is a surjection to the kernel of $\overline\tau:
E_0^{n,p}(B) \to\ E_0^{n,p}(C)$, and there is a split short exact
sequence $ 0 \to\ E_0^{n,p}(A)  \to\  E_0^{n,p}(B) \to\
E_0^{n,p}(C) \to\ 0 $. Statement (6) as well as the proposition
follows from the above.
\end{proof}

The next theorem follows at once as it was checked in \cite{C},
page $251$, that the local coefficient system for the extension

$$1 \to F_n \to P_{n+1}\to P_n \to 1$$ is trivial over the integers.

\begin{thm}\label{thm: restricted Lie algebras}
The restricted Lie algebra obtained from the mod-$p$ descending
central series for $P_{k}$ is given by $\mathcal{L}_{k,p} $  the
free restricted Lie algebra over $\mathbb Z/ p \mathbb Z$
generated by elements $B_{i,j}$ with $1 \leq i < j \leq k$, modulo
the infinitesimal braid relations:
\begin{description}
\item[(i)] $[B_{i,j}, B_{s,t}] =0 $ if $\{ i,j\} \cap \{ s,t\} =
\phi,$ \item[(ii)] $[ B_{i,j}, B_{i,t} + B_{t,j}] = 0$ if $1 \leq
i < t < j \leq k $, and \item[(iii)] $[B_{t,j}, B_{i,j} + B_{i,t}]
= 0$ if $1 \leq i < t < j \leq k$.
\end{description}

Furthermore there is a split short exact sequence of restricted
Lie algebras

\[
\begin{CD}
0 \to\ E_{0}^{*,p}(F_n) @>{ E_{0}^{*,p}(i) }>> E_{0}^{*,p}
(P_{n+1}) @>{E_{0}^{*,p}(d_n)}>> E_{0}^{*,p}(P_n)\to\ 0
\end{CD}
\] where $E_{0}^{*,p}(F_n)$ is the free restricted Lie algebra generated by
$B_{i, n+1}$ for $1 \leq i < n+1$. In addition,
$E_{0}^{*,p}(P_{n+1})$ is additively isomorphic to $
E_{0}^{*,p}(P_n) \oplus E_{0}^{*,p}(F_n).$
\end{thm}

The following is an ``integrality" statement concerning as
embeddings of Lie algebras may not induce an embedding after
mod-$p$ reduction.

\begin{thm}\label{thm: reduction mod p}
If $n \geq 1$, the induced maps $$E_0^*(\Theta_n): E_0^*(F_{n})
\to\ E_0^*(P_{n+1}),$$ and $$E_0^{*,p}(\Theta_n): E_0^{*,p}(F_{n})
\to\ E_0^{*,p}(P_{n+1})$$ are monomorphisms.
\end{thm}

\begin{proof}
Notice that the map $E_0^*(\Theta_n)$ sends a generator
$[\cdots[x_1,x_{j_1}]x_{j_2}]\cdots ] x_{j_{q-1}}]x_{j_q}]$ to
$$(\pm 1) [\cdots[\Lambda_n
,\gamma_{j_q}]\gamma_{j_{q-1}}]\cdots ]
\gamma_{j_{2}}]\gamma_{j_1}] + \Omega$$ where $\Omega$ projects to
an element of lower filtration degree in $H_1L[\Lambda_n^{C}]$ by
Theorem \ref{thm: filtration preserving on homology}.

By \ref{thm: isomorphism}, the map $$ E_0^*(\Theta_n):
L[x_1^{B}]\to\ L[\Lambda_n^{C}]$$ induces an isomorphism of Lie
algebras. Hence, the induced map on the level of universal
enveloping algebras is an isomorphism over the integers. Thus,
there is an induced isomorphism after reduction modulo $p$, and an
isomorphism on the level of restricted Lie algebras. This
suffices, and the Theorem follows.
\end{proof}

\section{On braid groups, and axioms for connected $CW$-complexes}\label{sec:axioms}

The purpose of this section is to give axioms which characterize
$CW$-complexes in terms of braid groups when viewed within the
context of simplicial groups. First consider the category of
groups $\mathcal G$, and the category of reduced simplicial groups
$\mathcal S \mathcal G$. Let $\mathcal C$ denote a small category.
The definition of a simplicial subgroup is used next, and is
defined below as the authors are unaware of an appropriate
reference.

Consider the following axioms: Let $\mathcal B$ denote the
smallest subcategory of $\mathcal S \mathcal G$ which satisfies
the following properties:

\begin{enumerate}
    \item The simplicial group $\mathrm{AP}_*$ is in $\mathcal B$.
    \item If $\Pi$, and $\Gamma$ are in $\mathcal B$, then the
    coproduct $\Pi \vee \Gamma$ is in $\mathcal B$.
    \item If $\Pi$ is in $\mathcal B$, and $\Gamma$ is a
    simplicial subgroup of $\Pi$, then $\Gamma$ is in $\mathcal B$.
     \item If $\Pi$ is in $\mathcal B$, and $\Gamma$ is a
     simplicial
    quotient of $\Pi$, then $\Gamma$ is in $\mathcal B$.
\end{enumerate}

The next result is stated in the Introduction as Theorem
\ref{thm:theorem one}.
\begin{thm}\label{thm:theorem CW axiomatization}
Let $X(i)$, $i = 1,2$ denote path-connected $CW$-complexes with a
continuous function $f:X(1) \to X(2)$. Then there exist elements
$\Gamma_{X(i)}$, together with a morphism $\gamma:\Gamma_{X(1)}
\to \Gamma_{X(2)}$ in $\mathcal B$ such that the loop space
$\Omega(X(i))$ is homotopy equivalent to the geometric realization
of $\Gamma_{X(i)}$, and the induced map $|\gamma|:|\Gamma_{X(1)}|
\to |\Gamma_{X(2)}|$ is homotopic to $\Omega(f)$.
\end{thm}

First, the definition of a simplicial subgroup is required. The
authors are unaware of a good reference, so additional features
are listed below.

\begin{defn}\label{defn: subgroup }
\begin{enumerate}
    \item A map $f:G \to H$ in $\mathcal C$ is a monomorphism provided
    whenever there are two maps $\alpha, \beta:\pi \to G$ in $\mathcal C$ such that
$$f\circ \alpha = f \circ \beta,$$ then
$$ \alpha = \beta.$$
    \item A map $f:G \to H$ in $\mathcal C$ is an injection
provided $f$ is one-to-one on the underlying sets.
\end{enumerate}

\end{defn} Next recall the following standard fact.

\begin{prop}\label{prop:injection}
A map $f:G \to H$ in $\mathcal G$ is an injection if and only if
$f$ is a monomorphism.
\end{prop}

\begin{proof}
Assume that $f$ is a monomorphism. It will be checked that
$\mathrm{ker}(f)$ is trivial. Assume that there is some
non-identity element $x$ in $\mathrm{ker}(f)$. Define

\begin{enumerate}
    \item $\alpha :\mathbb Z \to G$ by $\alpha(n) = 1$, and
    \item $\beta:\mathbb Z \to G$ by $\beta(n) = x^n$.
\end{enumerate} Then $f\circ \alpha = f \circ \beta$, but
$\alpha(x) \neq \beta(x)$, and so  $\alpha \neq \beta$. That is a
contradiction, and thus $\mathrm{ker}(f)$ is trivial, and so $f$
is an injection.

Next, assume that $f$ is an injection, and that $f\circ \alpha = f
\circ \beta$. Thus $f(\alpha(x)) = f (\beta(x))$. Hence $\alpha(x)
= \beta(x)$ and the proposition follows.
\end{proof}

The next step is to check an analogous statement for reduced
simplicial groups. First, two definitions should be given.

\begin{defn}\label{defn:simplicial subgroup }
Let $\Gamma$, and $\Pi$ denote reduced simplicial groups. Then
$\Gamma$ is a simplicial subgroup of $\Pi$ provided there is a
monomorphism $\phi: \Gamma \to \Pi$ in $\mathcal S \mathcal G$.
\end{defn}

\begin{defn}\label{defn:simplicial horseradish}
Let $\Gamma$, and $\Pi$ denote reduced simplicial groups. Then the
pair  $(\Pi,\Gamma)$ is a simplicial group pair provided there is
a morphism in $\mathcal S \mathcal G$ given by $\phi: \Gamma \to
\Pi$ which is a degree-wise injection of groups.
\end{defn}

\begin{prop}\label{prop:injections}
Given simplicial groups $\Gamma$, and $\Pi$, the following are
equivalent.
\begin{enumerate}
    \item $\Gamma$ is a simplicial subgroup of $\Pi$.
    \item The pair $(\Pi,\Gamma)$ is a simplicial group pair.
\end{enumerate}
\end{prop}

\begin{proof}
Assume that $f: \Gamma \to \Pi$ is a monomorphism in $\mathcal S
\mathcal G$. It will be checked that in each simplicial degree
$n$, $\mathrm{ker}(f)$ is trivial. Assume that $n$ is the minimal
degree for which there is a non-trivial element $x$ in
$\mathrm{ker}(f)$.

Let $$\bar G\langle x\rangle$$ denote the simplicial closure of
$x$ as given in \cite{CL}. Thus $\bar G\langle x\rangle$ is both a
simplicial group, and the natural morphism of simplicial groups
$\beta:\bar G\langle x\rangle \to \Gamma$ satisfies
$$\beta(x) = x.$$ Thus, there are morphisms in $\mathcal S
\mathcal G$ given by
\begin{enumerate}
    \item $\alpha :\bar G\langle x\rangle \to \Gamma$ by $\alpha(x) = 1$, and
    \item $\beta: \bar G\langle x\rangle \to \Gamma$ by $\beta(x) = x$.
\end{enumerate} Then $f\circ \alpha = f \circ \beta$, but
$\alpha(x) \neq \beta(x)$, and so  $\alpha \neq \beta$.

This statement contradicts the fact that $f$ is monomorphism. Thus
$\mathrm{ker}(f)$ is trivial, and $f$ is an injection in each
degree.

Next, assume that $f$ is an injection in each degree, and that
$f\circ \alpha = f \circ \beta$. Thus in each degree, it follows
that $f(\alpha(x)) = f (\beta(x))$. Hence $\alpha(x) = \beta(x)$
and the proposition follows.
\end{proof}

\section{Proof of Theorem \ref{thm:theorem one}}

The simplicial group $\mathrm{AP}_*$ is in $\mathcal B$ by axiom
$1$. By Theorem  \ref{thm:Embedding braids}, there is a morphism
of simplicial groups $$\Theta:F[S^1] \to \mathrm{AP}_*$$ which is
a degree-wise injection. Thus $F[S^1]$ is a simplicial subgroup of
$\mathrm{AP}_*$ by \ref{prop:injections}. Hence $F[S^1]$ is in
$\mathcal B$ by axiom $3$.

Notice that coproducts, $\mathrm{AP}_* \vee \mathrm{AP}_*$, as
well as $F[S^1]\vee F[S^1]$ are in $\mathcal B$ by axiom $2$.
Since $F[S^n]$ is a subgroup of $F[S^1]\vee F[S^1]$ for $n \geq
1$, $F[S^n]$ is in $\mathcal B$. By passage to coproducts, $\vee
_{n \epsilon T}F[S^n]$ is in $\mathcal B$ for any set $T$.

The next statement concerning push-outs is the simplicial analogue
of a classical result of J.~H.~C.~Whitehead \cite{Fie,Loday}.
found in \cite{WU}.

\begin{prop}\label{prop:simplicial push-outs}
Let $G_0$, $G_1$, and $G_2$ be simplicial groups in $\mathcal B$
together with morphisms $\alpha: G_0 \to G_1$, and  $\beta: G_0
\to G_2$ in $\mathcal B$. Then

\begin{enumerate}
    \item the push-out $\Pi$ of $\alpha: G_0 \to G_1$, and  $\beta:
G_0 \to G_2$ is in $\mathcal B$, and
    \item if both $\alpha$, and $\beta$ are monomorphisms in $\mathcal
    B$, the classifying space of $\Pi$
    is the push-out of the classifing space construction of
    $\alpha: G_0 \to G_1$, and  $\beta: G_0 \to G_2$.
\end{enumerate}
\end{prop}

Let $Y$ be the cofibre of a map $a: \vee_{n \epsilon T} S^{n} \to
X$, and assume that $\Omega X$ is homotopy equivalent to a
simplicial group $G$ in $\mathcal B$. Then consider the push-out
of groups

\[
\begin{CD}
F[\vee_{n \epsilon T} S^{n-1}] @>{f}>> G \\
 @V{g}VV         @VV{\Theta}V     \\
\vee _{S}  \mathrm{AP}_*    @>{}>> \Gamma\\
\end{CD}
\] where $f$ is the natural extension of the attaching map $a$, and
$g$ can be chosen to be a monomorphism. If $f$ is a monomorphism,
the push-out $\Gamma$ lies in the category $\mathcal B$, then the
geometric realization of $G$ is homotopy equivalent to $\Omega Y$
because $\Gamma$ is the free product with amalgamation by the
subgroup $F[\vee_S S^{n-1}]$, and the coproduct of $\vee_S
\mathrm{AP}_*$ is contractible.

It will be checked next that $f$ may be assumed to be a
monomorphism. First assume that $G'$ is homotopy equivalent to
$\Omega X$ and $$f': F[\vee_{n \epsilon T} S^{n-1}] \to G'$$
represents the looping of the attaching map. Observe that
$F[\vee_{n \epsilon T} S^{n-1}]$ embeds in $\vee_S \mathrm{AP}_*$
via the injection $$g: F[\vee_{n \epsilon T}  S^{n-1}]\to \vee_S
\mathrm{AP}_*.$$

Next, recall that $\mathcal B$ contains products as it contains
both coproducts, and quotients. Let $$G = G'\times
(\vee_S\mathrm{AP}_*),$$ and let $f$ be the injection $$f'\times
g: F[\vee_{n \epsilon T} S^{n-1}] \to G.$$ Notice that $G$ is
homotopy equivalent to $\Omega X$ because $\mathrm{AP}_*$ is
contractible, and that $G$ is in $\mathcal B$. Hence $f$ may be
replaced by $f'\times g$, and thus assumed to be a monomorphism.

To verify naturality, the previous arguments apply to continuous
maps $f:Y \to Z$ of connected $CW$ complexes with $Y = \vee_{n
\epsilon T}S^n$. In addition, if $Y$ be the cofibre of a map $a:
\vee_{n \epsilon T} S^{n} \to X$, an application of the push-out
constructed above suffices.

The result follows by \ref{prop:simplicial push-outs}.

\section{Appendix: a sample computation}
The purpose of this section is to list a sample computation for
the values of $\Theta_3$.

\begin{enumerate}
\item

$\Theta_3([x_1,x_2]x_2]x_3]) =
[\Lambda_3,\gamma_3]\gamma_2]\gamma_2]$,

\item $\Theta_3([x_1,x_2]x_3]x_2]) \cong
-[\Lambda_3,\gamma_2]\gamma_3]\gamma_2] +
2[\Lambda_3,\gamma_3]\gamma_2]\gamma_2]$ modulo decomposables, and

\item $\Theta_3([x_1,x_3]x_2]x_2]) \cong
[\Lambda_3,\gamma_2]\gamma_2]\gamma_3] -
2[\Lambda_3,\gamma_2]\gamma_3]\gamma_2] +
2[\Lambda_3,\gamma_3]\gamma_2]\gamma_2]$  modulo decomposables.
\end{enumerate}

\end{document}